# Width Laws and Spectral Geometry

*Universality, Geometric Memory, and Reconstruction*

Omri Abas

*Independent researcher, Kadima-Zoran, Israel*

omri.abas1@gmail.com

ORCID: 0009-0003-4857-0790

6 August 2026

---

**Abstract**

We develop a common framework for random width laws, spectral populations, and geometric reconstruction. For a $d$-dimensional orthotope, we prove an exact parity law for the maximal $\pi^{-1}$-grade of every spherical width cumulant, including noncancellation and sign in all dimensions and orders. The first $d$ scalar width moments recover the unordered side vector, and $d-1$ moments are generically insufficient. Each Laplace mode generates an auxiliary width law whose upper endpoint satisfies $M_{n,a} = \lambda_n(a)^{1/2}/\pi$. At high energy the modal coordinate partitions converge to a universal Dirichlet law, while an unsmoothed measure-valued cutoff expansion retains the first geometric memory at face scale. Its simplex moment determines, up to an explicit nonzero factor and a separate off-diagonal argument, a basis-independent projector-gradient Weyl tensor that reconstructs the orthotope. Genuine edge-scale jumps obstruct a third coefficient for the total raw cutoff; exact mixed-boundary Möbius inversion isolates every coordinate stratum and restores a recursive bulk–boundary expansion with a smaller remainder. Beyond orthotopes, we prove direction-labelled identifiability for a canonical linear–quadratic class and finite recovery from direction-sensitive ridge moments under a generator bound. In dimension three, a global great-circle incidence calculus gives the exact step, fold, endpoint-fold, and corner coefficients of reduced zonotopal width densities, including an explicit non-simple corner cancellation. The results distinguish universal aggregation, recoverable geometric memory, and the remaining scalar inverse problem.

**2020 Mathematics Subject Classification:** Primary 52A22; Secondary 35P20, 52B11, 52C35, 60D05.

**Keywords:** random widths; orthotopes; spectral geometry; geometric reconstruction; zonotopes; spherical pushforwards; lattice-point asymptotics.

---

## 1. Introduction and principal results

Let

$$K_a = \prod_{i=1}^{d} \left[-\frac{a_i}{2}, \frac{a_i}{2}\right], \qquad a = (a_1, \dots, a_d) \in (0,\infty)^d,$$

be an orthotope, and let

$$W_a(u) = \sum_{i=1}^{d} a_i\,|u_i|, \qquad u \in S^{d-1},$$

be its width in direction $u$. If $U$ is uniform on the sphere, then $W_a(U)$ is a scalar random variable obtained by forgetting the direction that produced the width. The Dirichlet and Neumann spectra of the same orthotope are generated by the reciprocal quadratic form

$$\lambda_n(a) = \pi^2 \sum_{i=1}^{d} \frac{n_i^2}{a_i^2},$$

with the usual positive or nonnegative index sets. These two objects appear different: one is a spherical pushforward law, the other a discrete spectral multiset. The purpose of this article is to identify the exact structure they share, to determine where geometric information survives aggregation, and to state sharply which forms of reconstruction are and are not justified.

The organizing principle is that geometric information occupies several levels. A direction-labelled width function may be replaced by its unlabelled law; an indexed family of auxiliary laws may be replaced by their endpoints; labelled endpoints may be replaced by an ordinary unlabelled spectrum; and a finite spectral population may be replaced by its high-energy limit. Every step forgets information. Nevertheless, the information is not simply lost. It may move from the universal bulk into a boundary-scale coefficient, into an isolated coordinate stratum, or into the incidence and singular coefficients of a spherical pushforward. The main results make this displacement explicit.

### *1.1. Width laws and the companion problem*

The distributional study of width has several distinct antecedents. Walters derived the projected-area distribution of a randomly oriented rectangular parallelepiped [29]. If its edge lengths are $b_1, b_2, b_3$, the projected area is

$$b_2 b_3 |u_1| + b_3 b_1 |u_2| + b_1 b_2 |u_3|.$$

The bijection

$$(b_1, b_2, b_3) \mapsto (a_1, a_2, a_3) = (b_2 b_3, b_3 b_1, b_1 b_2),$$

whose inverse is $b_1 = (a_2 a_3 / a_1)^{1/2}$ and cyclically, identifies this projected-area law with the absolute-linear-form width law of a rectangular box. Finch computed the mean and mean square width of the regular tetrahedron and the cube, noting that the density itself was not known [7]. Kabluchko, Litvak and Zaporozhets proved several regular-polytope moment formulas conjectured by Finch and established a Gaussian representation for higher random-width moments of arbitrary compact convex bodies [15]. Akiyama and Kamae developed a different but complementary planar direction, studying normalized width deviation and its extremizers for convex polygons [2].

The companion article derives the global width density of an arbitrary three-dimensional rectangular box and proves reconstruction from its first three cumulants [1]. The present article starts from that result but changes both dimension and question. It develops an all-dimensional moment and cumulant calculus, embeds a fixed orthotope into an infinite family of auxiliary width laws whose endpoints are spectral, and then follows geometric information through high-energy aggregation, boundary-scale correction, stratum isolation, and singular incidence. The direct derivations are retained even where a general representation theorem is available, because the coordinate factors and subset filtration are essential to the reconstruction and spectral arguments.

### *1.2. Principal results*

The first result is an exact arithmetic dichotomy for spherical width cumulants. If $\kappa_N^{(d)}(a)$ denotes the $N$-th cumulant of $W_a(U)$, Theorem 3.1 proves, for every positive side vector and every $N \geq 1$, that

$$\deg_{\pi^{-1}} \kappa_N^{(d)}(a) = \begin{cases} N, & d \text{ even}, \\ \lfloor N/2 \rfloor, & d \text{ odd}, \end{cases}$$

and determines the sign of every extremal coefficient. The odd-dimensional, odd-order case is the nontrivial part: a symmetric half-plane polarization, a Pochhammer multiplier, and a positive Stieltjes–beta transfer rule out cancellation for every order and every positive side vector. Theorem 3.5 then shows that the first $d$ scalar width moments determine the unordered side vector of a $d$-orthotope, while Proposition 3.6 proves that $d-1$ moments are generically insufficient. Thus finite moment closure is both exact and generically sharp within the orthotopal class.

Corollary 3.4 extracts more than the degree statement: it gives the extremal coefficients explicitly at both parities and proves that their normalized root-exponential growth is governed by the smallest positive reciprocal root of the transformed even polynomial $P_a$. In dimension three this collapses to the single-root formulas underlying the companion result.

The second result is the spectral-family construction. For a nonzero mode index $n$, define

$$W_{n,a}(u) = \sum_{i=1}^{d} \frac{n_i}{a_i} |u_i|, \qquad M_{n,a} = \max_{u \in S^{d-1}} W_{n,a}(u).$$

Proposition 4.1 gives the endpoint identity

$$\lambda_n(a) = \pi^2 M_{n,a}^2.$$

Consequently, the Dirichlet or Neumann spectrum is the multiset of the quantities $\pi^2 M_{n,a}^2$, the squared upper endpoints of a mode-indexed family of auxiliary width laws multiplied by $\pi^2$. This identity is exact, but it is also an information boundary: an endpoint does not contain the interior of its law, and an ordinary spectrum additionally forgets the mode labels. Section 4 separates indexed laws, labelled endpoints, and unlabelled spectra before giving finite reconstruction at each justified level. In particular, the $d$ nonconstant algebraic coefficients of the complete Dirichlet heat trace recover the elementary symmetric functions of the side lengths and hence the orthotope.

The third result resolves the apparent conflict between spectral universality and geometric memory. Theorem 5.2 proves that the high-energy distribution of the modal coordinate partitions converges to the universal Dirichlet$(1/2, \dots, 1/2)$ law, independently of the aspect ratio. This uses a classical directional lattice-point limit, but identifies its pushforward with the same spherical-square law that controls the continuous width problem. Corollary 5.3 strengthens the connection: the empirical distribution of the corresponding labelled maximizing directions, evaluated through the original side vector, converges to the width law $\mathcal{L}(W_a(U))$. The correspondence is therefore exact at the level of labelled modal directions, not at the level of an ordinary unlabelled spectrum. Theorem 5.1 also separates fixed-label population optimization from ordered eigenvalue optimization: no orthotope minimizes every prescribed shape-independent labelled population, whereas every orthotope is uniquely optimal for a suitable finite population. The cube reappears in Theorem 5.4 as the unique orthotope with exact coordinate balance at every cutoff and every heat scale.

Section 6 gives the central universality–memory theorem. Theorem 6.2 obtains an unsmoothed, measure-valued sharp-cutoff expansion consisting of a universal bulk law, an exact face-scale correction, and a remainder smaller than face scale. The singular part of the correction is supported on the simplex faces. In the principal frame its coordinate moments satisfy

$$\int_{\Delta_{d-1}} q_i \, d\mathfrak{M}_{a,s}^{\mathrm{sing}} = -\frac{\gamma_{d,s}}{2(d-1)V(a)} \, \mathsf{F}(K_a)_{ii}.$$

The off-diagonal entries vanish by the separate sine–cosine argument of Subsection 6.2. Thus these moments determine, up to the displayed explicit nonzero factor, the basis-independent projector-gradient tensor

$$\mathsf{F}(K_a) = \int_{\partial K_a} (I - \nu \otimes \nu) \, dS,$$

which explicitly reconstructs every side length. Theorem 6.4 then proves that the next term of the total raw cutoff is obstructed by genuine edge-scale jumps; the absence of a deterministic coefficient is not a technical defect of the estimate. Exact mixed Dirichlet–Neumann Möbius transforms resolve the obstruction. Theorem 6.5 isolates every coordinate stratum, and Theorem 6.6 shows that within an isolated stratum of intrinsic dimension at least two the bulk–boundary expansion restarts in that dimension, with its own next coefficient and a strictly smaller remainder. At edge level this two-term restart begins in ambient dimension $d \geq 4$; when $d = 3$, the isolated one-dimensional edge stratum has only the leading term and the remainder stated in Theorem 6.6. The hierarchy is therefore recursive:

universal bulk → face memory → raw edge obstruction → exact stratum isolation → nested boundary memory.

This result is stronger than a scalar Weyl coefficient: it identifies the directional coefficient, converts it into a geometric tensor, proves explicit reconstruction, and locates the precise scale at which the uncancelled sharp cutoff ceases to possess a deterministic next term. The mixed-boundary raw data include the pure Dirichlet spectrum. The finite heat-invariant system of Section 4 and the finite stratum-coefficient system of Section 6 recover the same side geometry through different extraction mechanisms.

The fourth result establishes an exact identifiability boundary beyond orthotopes. For the canonical linear–quadratic class

$$\mathcal{W}_{G,A}(u) = \sum_j |g_j \cdot u| + 2\sqrt{u^T A u},$$

Theorem 7.1 proves that the complete direction-labelled function determines the reduced generator family and the positive-semidefinite quadratic form, up to permutation and sign changes of the generators; equivalently, it determines the centred body, and without a prescribed centre determines the body up to translation. Under a fixed generator bound, Theorem 7.3 replaces the continuum by finitely many direction-sensitive ridge moments and one quadratic residual tensor moment. These inputs are not ordinary scalar moments of $\mathcal{W}_{G,A}(U)$, and they are not ordinary spectral aggregates.

The distinction is essential: Subsection 7.1 gives a noninjectivity result for unlabelled scalar width laws on an unrestricted smooth centrally symmetric class. Within distinguished subclasses, however, less data suffice: the complete scalar width law reconstructs a right circular cylinder, and finitely many even scalar moments reconstruct an ellipsoid.

The fifth result is a global singularity calculus for reduced zonotopal width laws in dimension three. Section 8 constructs the density by co-area on the complete great-circle arrangement and derives the exact incidence conditions and local coefficients for steps, folds, endpoint folds, and corners. Theorem 8.1 distinguishes what is unconditional from what is generic: steps are sign-rigid, isolated folds and simple corners are nonzero, and resonant contributions superpose algebraically. An explicit non-simple arrangement has six nonzero local slopes but zero total corner coefficient. This example corrects the universal corner clause conjectured in the companion article. The replacement theorem is stronger in the appropriate sense: it gives the exact cyclic coefficient, characterizes the generic noncancellation regime, and exhibits the mechanism by which the universal corner assertion fails at a resonant non-simple corner.

### *1.3. Scope and information boundaries*

The reconstruction theorems in this article are deliberately data-specific. For orthotopes, the first $d$ scalar width moments and the $d$ nonconstant algebraic heat invariants determine the unordered side multiset, whereas the face-memory tensor and the isolated codimension-two coefficients determine the coordinate-labelled side vector, together with its principal-axis frame when the side lengths are distinct. These are distinct observables, although their information content becomes equivalent after the corresponding reconstruction maps are applied. For the mixed linear–quadratic class, complete direction-labelled data, or the finite direction-sensitive data of Theorem 7.3 under a generator bound, determine the canonical geometry. No theorem here asserts that finitely many ordinary scalar moments identify the full mixed class, that an ordinary unlabelled spectrum contains the auxiliary-law interiors, or that a scalar width law determines an unrestricted convex body.

Conjecture 8.2 isolates the remaining finite scalar problem. It asks whether, on each bounded canonical linear–quadratic class, a dimension- and complexity-dependent finite number of ordinary scalar moments is generically injective outside a proper exceptional set. The conjecture is consistent with the finite subclass theorems and with the direction-sensitive reconstruction, but neither implies it. It is retained as an explicit frontier rather than used as part of the proved framework.

### *1.4. Organization*

Section 2 develops the common spherical and Boolean-lattice calculus. Section 3 proves cumulant parity and finite orthotope reconstruction. Section 4 introduces the auxiliary width family and separates its spectral data levels. Section 5 proves partition universality, modal-width correspondence, conditional optimization, and all-scale balance. Section 6 derives the sharp bulk–face expansion, face-memory tensor, raw edge obstruction, and exact mixed-boundary stratum tomography. Section 7 establishes the positive and negative identifiability results beyond orthotopes. Section 8 gives the stratified zonotopal density theorem, the cancellation counterexample, and the finite scalar conjecture. Section 9 assembles the consequences and records the proved information boundaries.

# 2. Spherical widths and Boolean-lattice decompositions

Fix an integer $d \geq 2$, and write $[d] = \{1, \dots, d\}$. Let

$$K_a = \prod_{i=1}^{d} \left[-\frac{a_i}{2}, \frac{a_i}{2}\right], \qquad a = (a_1, \dots, a_d) \in (0, \infty)^d,$$

be the centred orthotope with side vector $a$. Its width in the deterministic direction $u \in S^{d-1}$ is

$$W_a(u) = \sum_{i=1}^{d} a_i \, |u_i|. \qquad (2.1)$$

Throughout, $u$ denotes a deterministic unit vector, whereas $U$ denotes a random vector distributed uniformly on $S^{d-1}$. The random width is therefore $W_a(U)$.

The purpose of this section is twofold. First, we establish the exact spherical moment calculus used in the reconstruction theorem of Section 3. Second, we isolate two decompositions indexed by the Boolean lattice $2^{[d]}$: a continuous decomposition according to the support of odd exponents, and a discrete decomposition according to the coordinates constrained to vanish. Their common indexing is structural, not analogical. Lemma 2.4 will later be applied directly to absolutely convergent spectral sums, thereby resolving them into contributions from coordinate strata.

### *2.1. Squared spherical coordinates*

For $U = (U_1, \dots, U_d)$ uniform on $S^{d-1}$, set

$$X_i = U_i^2, \qquad 1 \leq i \leq d.$$

The vector $X = (X_1, \dots, X_d)$ takes values in the closed simplex

$$\overline{\Delta}_{d-1} = \{x \in [0, \infty)^d : \sum_{i=1}^{d} x_i = 1\}.$$

Its density is supported on the relative interior

$$\Delta_{d-1} = \{x \in (0, \infty)^d : \sum_{i=1}^{d} x_i = 1\}.$$

**Proposition 2.1.** The squared-coordinate vector $X$ has the Dirichlet distribution with parameter vector $(1/2, \dots, 1/2)$. Thus its density on $\Delta_{d-1}$, with respect to the coordinate measure $dx_1 \cdots dx_{d-1}$, where

$$x_d = 1 - \sum_{i=1}^{d-1} x_i,$$

is

$$f_X(x) = \frac{\Gamma(d/2)}{\Gamma(1/2)^d} \prod_{i=1}^{d} x_i^{-1/2}.$$

**Proof.**
Let $G_1, \dots, G_d$ be independent standard Gaussian variables and put

$$R = \left( \sum_{i=1}^{d} G_i^2 \right)^{1/2}.$$

Rotational invariance gives

$$U \overset{d}{=} \frac{G}{R}.$$

The variables $Y_i = G_i^2/2$ are independent gamma variables with common shape parameter $1/2$ and common scale 1. Writing

$$T = \sum_{i=1}^{d} Y_i, \qquad X_i = \frac{Y_i}{T},$$

and applying the change of variables

$$(y_1, \dots, y_d) = (tx_1, \dots, tx_d),$$

whose Jacobian on the simplex coordinates is $t^{d-1}$, separates the joint density into a gamma density in $t$ and the displayed Dirichlet density in $x$. Since $X_i = G_i^2/R^2 = U_i^2$, the assertion follows. □

### *2.2. Spherical and radial moments*

Write

$$m_N(a) = \mathbb{E}\, W_a(U)^N, \qquad N \geq 0.$$

The following lemma gives both the radial factorization and the exact multinomial expression for $m_N(a)$.

**Lemma 2.2.** Let $G = (G_1, \dots, G_d)$ be standard Gaussian in $\mathbb{R}^d$, let $R = \| G \|_2$, and let

$$S_a = \sum_{i=1}^{d} a_i \, |G_i|.$$

Then $R$ and $U = G/R$ are independent,

$$S_a = R\, W_a(U),$$

and, for every integer $N \geq 0$,

$$\mathbb{E}S_a^N = \mathbb{E}R^N\, m_N(a), \qquad \mathbb{E}R^N = 2^{N/2}\frac{\Gamma\big((d+N)/2\big)}{\Gamma(d/2)}.$$

Consequently,

$$m_N(a) = \frac{\Gamma(d/2)\,N!}{\Gamma\big((d+N)/2\big)} \sum_{\substack{q_1,\dots,q_d\geq 0\\ q_1+\cdots+q_d=N}} \prod_{i=1}^{d} \frac{\Gamma\big((q_i+1)/2\big)}{\Gamma(1/2)\,q_i!}\, a_i^{q_i}. \qquad (2.2)$$

**Proof.**
The Gaussian polar decomposition gives $G = RU$, where $R$ and $U$ are independent and $U$ is uniform on $S^{d-1}$. Hence

$$S_a = \sum_{i=1}^{d} a_i\,|G_i| = R\sum_{i=1}^{d} a_i\,|U_i| = R\,W_a(U).$$

The radial density is proportional to $r^{d-1}e^{-r^2/2}$. Direct integration yields

$$\mathbb{E}R^N = \frac{\int_0^\infty r^{d+N-1}\,e^{-r^2/2}\,dr}{\int_0^\infty r^{d-1}\,e^{-r^2/2}\,dr} = 2^{N/2}\frac{\Gamma\big((d+N)/2\big)}{\Gamma(d/2)}.$$

Independence now gives the first moment identity.

Alternatively, expanding (2.1) multinomially and applying Proposition 2.1 gives

$$m_N(a) = \sum_{|q|=N}\binom{N}{q_1,\dots,q_d} a^q\,\mathbb{E}\prod_{i=1}^{d} X_i^{q_i/2}.$$

The Dirichlet moment formula

$$\mathbb{E}\prod_{i=1}^{d} X_i^{q_i/2} = \frac{\Gamma(d/2)}{\Gamma\big((d+N)/2\big)}\prod_{i=1}^{d}\frac{\Gamma\big((q_i+1)/2\big)}{\Gamma(1/2)}$$

then proves (2.2). □

The Gaussian radial identity is the orthotopal specialization of the general higher-moment Gaussian-width formula for compact convex bodies established by Kabluchko, Litvak and Zaporozhets [15]. We retain the direct derivation because the proof below requires the individual coordinate factors visible in (2.2).

### *2.3. The continuous subset decomposition*

For a multi-index $q = (q_1,\dots,q_d)$, define its odd support by

$$\operatorname{odd}(q) = \{i \in [d] : q_i \text{ is odd}\}.$$

Because $|q| = N$,

$$|\operatorname{odd}(q)| \equiv N \pmod 2.$$

Moreover,

$$\frac{\Gamma\big((q_i+1)/2\big)}{\Gamma(1/2)}$$

is rational when $q_i$ is even and is a rational multiple of $\pi^{-1/2}$ when $q_i$ is odd. Since the cardinality of $\mathrm{odd}(q)$ has the same parity as $N$, the total power of $\pi^{-1}$ is integral.

**Lemma 2.3.** For every $N \geq 0$, there exist homogeneous polynomials $M_{N,J}(a)$ with rational coefficients, indexed by subsets $J \subseteq [d]$ satisfying

$$|J| \leq N, \qquad |J| \equiv N \ (\mathrm{mod}\ 2),$$

such that

$$m_N(a) = \sum_{\substack{J\subseteq[d]\\ |J|\equiv N\ (\mathrm{mod}\ 2)}} \pi^{-\nu_{d,N}(|J|)}\, M_{N,J}(a). \qquad (2.3)$$

Here

$$\nu_{d,N}(s) = \begin{cases} s/2, & N \text{ even},\\ (s-1)/2, & N \text{ odd and } d \text{ odd},\\ (s+1)/2, & N \text{ odd and } d \text{ even}.\end{cases}$$

The polynomial $M_{N,J}$ is obtained by summing precisely those terms of (2.2) whose odd support is $J$.

**Proof.**
Partition the multi-indices $q$ occurring in (2.2) according to $J = \mathrm{odd}(q)$. Every coordinate in $J$ contributes at least one to $|q| = N$, so $|J| \leq N$; the parity condition was proved above. The coordinate factors contribute $\pi^{-|J|/2}$. If $N$ is even, the outer gamma ratio in (2.2) is rational. If $N$ is odd, it is a rational multiple of $\sqrt{\pi}$ when $d$ is odd and of $\pi^{-1/2}$ when $d$ is even. These three cases give $\nu_{d,N}$. All remaining factors are rational factorial expressions, and homogeneity of degree $N$ follows from $|q| = N$. □

Thus the maximal possible power of $\pi^{-1}$ is governed by the largest admissible odd support. This filtration is the starting point of the cumulant-parity theorem.

### *2.4. The discrete subset decomposition*

Let $f \in \ell^1(\mathbb{Z}^d)$ be coordinatewise even. For $J \subseteq [d]$, define

$$L_J(f) = \sum_{m\in\mathbb{Z}^J} f\,(m_J, 0_{J^c}).$$

The next identity is stated for absolutely summable functions, rather than only finitely supported ones, because Subsection 6.6 applies it to heat-weighted spectral sums.

**Lemma 2.4.** With $\mathbb{N} = \{1,2,\dots\}$,

$$\sum_{n\in\mathbb{N}^d} f\,(n) = 2^{-d} \sum_{J\subseteq[d]} (-1)^{d-|J|}\, L_J(f). \qquad (2.4)$$

**Proof.**
If $h \in \ell^1(\mathbb{Z})$ is even, then

$$\sum_{n\geq 1} h\,(n) = \frac{1}{2}\left(\sum_{m\in\mathbb{Z}} h\,(m) - h(0)\right).$$

Apply this identity successively in the $d$ coordinates. Choosing the full lattice sum in precisely the coordinates of $J$ and the zero-coordinate subtraction in the remaining $d - |J|$ coordinates gives $L_J(f)$ with sign $(-1)^{d-|J|}$. Absolute convergence justifies every use of Fubini and every rearrangement, proving (2.4). □

Equations (2.3) and (2.4) describe different objects: the first decomposes continuous spherical moments by odd coordinate support; the second decomposes a positive-orthant lattice sum by zero-coordinate strata. The link becomes operative only when both are applied to the same indexed family of orthotopal width laws. Section 3 develops the continuous reconstruction mechanism. Section 4 introduces that family, Section 5 applies (2.4) to finite cutoff populations, and Section 6 applies it both to raw cutoffs and to absolutely convergent heat populations while retaining the lower-dimensional coordinate strata.

# 3. Cumulant parity and finite reconstruction

Let

$$\kappa_N^{(d)}(a) = \kappa_N\big(W_a(U)\big)$$

denote the $N$-th cumulant of the random width. All coefficient extraction in this section is formal. Throughout the construction, $a_1, \dots, a_d$ are treated as algebraically independent indeterminates. Introduce a fresh indeterminate $\Xi$. Starting from the explicit gamma-function expressions above and applying the cumulant recursion in the polynomial ring, replace each displayed factor $\pi^{-j}$ by $\Xi^j$, before any side vector is evaluated. This defines canonical formal lifts

$$m_N(a),\ \kappa_N^{(d)}(a) \in \mathbb{Q}[a_1, \dots, a_d][\Xi].$$

We continue to write $\deg_{\pi^{-1}}$ and $[\pi^{-j}]$ for degree and coefficient extraction in $\Xi$. Both operations refer to these canonical formal lifts before evaluation of the side vector. When a coefficient is subsequently evaluated at $a \in (0,\infty)^d$, it is the corresponding coefficient polynomial in $\mathbb{Q}[a_1, \dots, a_d]$ that is evaluated. The first question concerns the exact maximal grade, not only its occurrence.

The main result of this section is stated before the machinery used to prove it.

**Theorem 3.1 (cumulant-parity theorem).** For every $a \in (0,\infty)^d$ and every $N \geq 1$,

$$\deg_{\pi^{-1}} \kappa_N^{(d)}(a) = \begin{cases} N, & d \text{ even,} \\ \lfloor N/2 \rfloor, & d \text{ odd.} \end{cases} \qquad (3.1)$$

The extremal coefficient never vanishes. More precisely:

1. if $d$ is even, then

$$\operatorname{sgn}\left([\pi^{-N}]\kappa_N^{(d)}(a)\right) = (-1)^{N-1};$$

2. if $d = 2r + 1$, then, for $M \geq 1$,

$$\operatorname{sgn}\left([\pi^{-M}]\kappa_{2M}^{(d)}(a)\right) = (-1)^{M-1},$$

and, for $M \geq 0$,

$$\operatorname{sgn}\left([\pi^{-M}]\kappa_{2M+1}^{(d)}(a)\right) = (-1)^{M}.$$

The proof occupies Subsections 3.1–3.3. The even-dimensional case is elementary. In odd dimension, even cumulants are controlled by a real-rooted polynomial, whereas odd cumulants require a positive Stieltjes representation and a half-step multiplier transfer.

### *3.1. The highest $\pi^{-1}$-grade*

Define the exponential moment series

$$\mathcal{M}_a(t) = 1 + \sum_{n\geq 1} m_n(a)\frac{t^n}{n!}.$$

By the moment–cumulant relation,

$$\log\mathcal{M}_a(t) = \sum_{n\geq 1} \kappa_n^{(d)}(a)\frac{t^n}{n!}.$$

Suppose first that $d = 2r + 1$. Introduce

$$P_a(z) = \sum_{k=0}^{r} \frac{\Gamma(d/2)}{\Gamma(d/2+k)} e_{2k}(a)\, z^k$$

and

$$R_a(z) = \sum_{k=0}^{r} \frac{\Gamma(d/2)}{\sqrt{\pi}\,\Gamma(d/2+k+1/2)} e_{2k+1}(a)\, z^k,$$

where $e_j(a)$ denotes the $j$-th elementary symmetric polynomial in $a_1, \dots, a_d$.

To make the filtration explicit, write

$$\frac{m_n(a)}{n!} = \sum_{j\geq 0} c_{n,j}(a)\,\pi^{-j}$$

and introduce the marked series

$$\widehat{\mathcal{M}}_a(t,y) = 1 + \sum_{n\geq 1}\sum_{j\geq 0} c_{n,j}(a) t^n y^j.$$

Assign the monomial $t^n y^j$ the defect $n - 2j$. By Lemma 2.3, defects are nonnegative. The terms of defect 0 arise exactly from even moments with $2k$ distinct coordinates carrying exponent 1; the terms of defect 1 arise exactly from odd moments with $2k + 1$ such coordinates. Hence

$$\widehat{\mathcal{M}}_a(t,y) = P_a(yt^2) + tR_a(yt^2) + \{\text{terms of defect at least 2}\}.$$

Taking the formal logarithm preserves the defect filtration. A term of defect 0 can use only $P_a(yt^2)$, while a term of defect 1 must contain exactly one factor from $tR_a(yt^2)$ and otherwise only defect-zero factors. Thus

$$\frac{[\pi^{-N}]\kappa_{2N}^{(2r+1)}(a)}{(2N)!} = [z^N]\log P_a(z), \qquad \frac{[\pi^{-N}]\kappa_{2N+1}^{(2r+1)}(a)}{(2N+1)!} = [z^N]\frac{R_a(z)}{P_a(z)}. \tag{3.2}$$

Products containing three or more odd-moment factors have defect at least 3 and therefore cannot enter the second identity.

The fact that $\deg P_a = r$ does not restrict (3.2) to $N \le r$. Since $P_a(0) = 1$, the formal series $\log P_a$ and $P_a^{-1}$ generally have coefficients of every order. For $N > r$, the extremal coefficient is produced by repeated products of lower-degree terms, equivalently by set partitions containing repeated pair blocks. This is precisely why the logarithmic formulation is required.

Suppose now that $d = 2r$. Lemma 2.3 gives

$$\deg_{\pi^{-1}} m_1 = 1, \qquad \deg_{\pi^{-1}} m_s < s \quad (s \ge 2).$$

For a partition $\mathcal{P}$ of $[N]$, the corresponding product in the moment–cumulant formula therefore satisfies

$$\deg_{\pi^{-1}} \prod_{B\in\mathcal{P}} m_{|B|} < \sum_{B\in\mathcal{P}} |B| = N$$

whenever $\mathcal{P}$ contains a nonsingleton block. Thus degree $N$ can occur only for the all-singleton partition, independently of whether $N \le d$. Put

$$c_{2r} = \frac{4^r r!\,(r-1)!}{(2r)!} > 0.$$

Then $m_1(a) = c_{2r}\pi^{-1}e_1(a)$. In the moment–cumulant formula, the power $\pi^{-N}$ can arise only from the partition of $[N]$ into $N$ singletons. Its coefficient is therefore

$$[\pi^{-N}]\kappa_N^{(2r)}(a) = (-1)^{N-1}(N-1)!\,(c_{2r}e_1(a))^N.$$

This proves the even-dimensional part of Theorem 3.1, including its sign assertion.

### *3.2. Real-rootedness and the even cumulants*

We now fix $d = 2r + 1$. Define

$$E_a(z) = \sum_{k=0}^{r} e_{2k}(a)z^k, \qquad O_a(z) = \sum_{k=0}^{r} e_{2k+1}(a)z^k.$$

The factorization

$$\prod_{i=1}^{2r+1}(1 + a_i t) = E_a(t^2) + tO_a(t^2)$$

controls their zeros.

For $y > 0$, write

$$\prod_{i=1}^{2r+1}(1+ia_iy)=A(y)e^{i\phi(y)}, \qquad \phi(y)=\sum_{i=1}^{2r+1}\arctan(a_iy),$$

where $A(y)>0$. Then

$$E_a(-y^2)=A(y)\cos\phi(y), \qquad yO_a(-y^2)=A(y)\sin\phi(y).$$

The function $\phi$ is strictly increasing,

$$\phi(0)=0, \qquad \phi(y)\uparrow\frac{(2r+1)\pi}{2} \quad \text{as } y\to\infty.$$

Consequently, the equations

$$\phi(y)=\left(j-\frac{1}{2}\right)\pi, \qquad 1\le j\le r,$$

give exactly $r$ finite positive zeros of $E_a(-y^2)$, while

$$\phi(y)=j\pi, \qquad 1\le j\le r,$$

give exactly $r$ finite positive zeros of $O_a(-y^2)$. Since both polynomials have degree $r$, all their zeros have thereby been accounted for. Their zeros are simple, negative and strictly interlacing.

The polynomial $P_a$ is obtained from $E_a$ by the coefficient multiplier

$$z^k\mapsto\frac{\Gamma(d/2)}{\Gamma(d/2+k)}\,z^k.$$

The following first-principles polarization lemma supplies the required preservation mechanism. It is a symmetric multiaffine form of the classical Grace–Walsh–Szegő principle; see [5, 12, 24, 28] for its classical setting.

**Lemma 3.2 (symmetric half-plane polarization).** Let

$$f(z)=\sum_{k=0}^{n}c_k\,z^k$$

be a polynomial of degree $n$ having no zero in

$$\mathbb{H}=\{z\in\mathbb{C}:\mathrm{Im}z>0\}.$$

Then its symmetric multiaffine polarization

$$\mathrm{Pol}_nf(z_1,\dots,z_n)=\sum_{k=0}^{n}\frac{c_k}{\binom{n}{k}}e_k(z_1,\dots,z_n)$$

does not vanish on $\mathbb{H}^n$.

**Proof.**
For a polynomial $p$ of exact degree $m$, define its polar derivative with pole $\alpha$ by

$$D_\alpha^{(m)}p(z)=mp(z)+(\alpha-z)p'(z).$$

Assume that the zeros of $p$ lie in the closed lower half-plane and that $\alpha \in \mathbb{H}$. Set

$$z = \alpha + \frac{1}{w}, \qquad q(w) = w^m p\left(\alpha + \frac{1}{w}\right).$$

Writing $\alpha = A + iB$, with $B > 0$, shows that the real boundary in the $z$-plane is mapped to

$$u^2 + \left(v - \frac{1}{2B}\right)^2 = \frac{1}{4B^2}, \qquad w = u + iv.$$

The closed lower half-plane is therefore mapped to the corresponding closed disk. Moreover,

$$q'(w) = w^{m-1} D_\alpha^{(m)} p\left(\alpha + \frac{1}{w}\right).$$

Every zero of $q'$ lies in the convex hull of the zeros of $q$. Indeed, if $w_0$ lay outside that convex hull, a separating line and a rotation would make all numbers $\left(w_0 - \zeta_j\right)^{-1}$, with $\zeta_j$ a zero of $q$, have positive real part. Hence

$$\frac{q'(w_0)}{q(w_0)} = \sum_j \frac{1}{w_0 - \zeta_j} \neq 0.$$

Thus the zeros of $q'$ remain in the disk, and the zeros of $D_\alpha^{(m)} p$ remain in the closed lower half-plane. The polar derivative has exact degree $m - 1$: if

$$p(z) = a_m \prod_{j=1}^{m} (z - \zeta_j),$$

then the coefficient of $z^{m-1}$ in $D_\alpha^{(m)} p$ is

$$a_m \left( m\alpha - \sum_{j=1}^{m} \zeta_j \right),$$

which cannot vanish because $m\alpha \in \mathbb{H}$, whereas $\sum_{j=1}^{m} \zeta_j$ lies in the closed lower half-plane.

Apply successively the polar derivatives with poles $z_1, \dots, z_n$. A direct induction on the monomial $z^k$ gives

$$D_{z_n}^{(1)} D_{z_{n-1}}^{(2)} \cdots D_{z_1}^{(n)} z^k = k!\,(n-k)!\, e_k(z_1, \dots, z_n).$$

Therefore

$$D_{z_n}^{(1)} D_{z_{n-1}}^{(2)} \cdots D_{z_1}^{(n)} f = n!\, \mathrm{Pol}_n f(z_1, \dots, z_n).$$

Each polar derivative preserves exclusion of $\mathbb{H}$, so the left-hand side cannot vanish when every $z_j \in \mathbb{H}$. □

We also require the following Pochhammer multiplier. It is a classical multiplier-sequence instance of Pólya–Schur theory [4]; the proof is included to keep the argument self-contained.

**Lemma 3.3 (Pochhammer multiplier).** Let $b > 0$, and let

$$p(z) = \sum_{k=0}^{n} c_k \, z^k$$

be a polynomial of exact degree $n$ having no zero in $\mathbb{H}$. Define

$$T_b p(z) = \sum_{k=0}^{n} \frac{c_k}{(b)_k} z^k, \qquad (b)_k = \frac{\Gamma(b+k)}{\Gamma(b)}. \qquad (3.3)$$

Then $T_b p$ has no zero in $\mathbb{H}$. Consequently, if $p$ has real coefficients and only real zeros, then $T_b p$ also has only real zeros. If all zeros of $p$ are negative and its coefficients are positive, the same is true of $T_b p$.

**Proof.**
Consider

$$J_{n,b}(z) = \sum_{k=0}^{n} \binom{n}{k} \frac{z^k}{(b)_k} = \frac{n!}{(b)_n} L_n^{b-1}(-z).$$

Rodrigues' identity gives

$$L_n^{b-1}(x) = \frac{x^{1-b} e^x}{n!} \frac{d^n}{dx^n} \left(e^{-x} x^{n+b-1}\right).$$

If $q$ is a polynomial of degree less than $n$, integration by parts $n$ times gives

$$\int_0^\infty L_n^{b-1}(x) q(x) x^{b-1} e^{-x} \, dx = 0.$$

All boundary terms vanish because $b > 0$. If $L_n^{b-1}$ had fewer than $n$ sign changes on $(0, \infty)$, choosing $q$ as the product of its sign-change points, with the appropriate overall sign, would make the integrand nonnegative and not identically zero. This contradicts the displayed orthogonality. Thus $L_n^{b-1}$ has $n$ simple positive zeros, and

$$J_{n,b}(z) = \prod_{j=1}^{n} \left(1 + \beta_j z\right), \qquad \beta_j > 0.$$

Let $F = \mathrm{Pol}_n p$. Comparison of coefficients gives

$$e_k(\beta_1, \dots, \beta_n) = \binom{n}{k} \frac{1}{(b)_k}.$$

Consequently, the coefficientwise map (3.3) satisfies

$$F(\beta_1 z, \dots, \beta_n z) = T_b p(z).$$

If $z \in \mathbb{H}$, then every $\beta_j z \in \mathbb{H}$, and Lemma 3.2 shows that $T_b p(z) \neq 0$. This proves preservation of upper-half-plane stability. If $p$ has real coefficients and only real zeros, then $T_b p$ has real coefficients. Any nonreal zero in the lower half-plane would have a conjugate zero in $\mathbb{H}$, so all zeros are real. If the coefficients of $p$ are positive, so are those of $T_b p$, and nonnegative roots are excluded. □

Since

$$P_a = T_{d/2} E_a,$$

Lemmas 3.2 and 3.3 imply that $P_a$ has only negative roots. Write

$$P_a(z) = \prod_{\ell=1}^{r}\left(1 + \frac{z}{\rho_\ell}\right), \qquad \rho_\ell > 0,$$

with multiplicities included. Then

$$[z^N]\log P_a(z) = \frac{(-1)^{N-1}}{N}\sum_{\ell=1}^{r} \rho_\ell^{-N}.$$

The coefficient is nonzero and has sign $(-1)^{N-1}$ for every $N \geq 1$, including $N > r$. By the first identity in (3.2), this proves the odd-dimensional, even-order part of Theorem 3.1.

### *3.3. Odd cumulants: the Stieltjes–beta mechanism*

Put

$$b = r + \frac{1}{2}, \qquad c_r = \frac{\Gamma(r + 1/2)}{\sqrt{\pi}\,\Gamma(r + 1)},$$

so that

$$P_a = T_b E_a, \qquad R_a = c_r T_{b+1/2} O_a.$$

Introduce the common-multiplier polynomial

$$B_a = T_b O_a.$$

We first derive a Stieltjes representation for $B_a/P_a$. The general relation between interlacing and rational $R$-functions is discussed in [25]; the signs required here are obtained directly. If $-s_j^2$ is a zero of $E_a$, then differentiation of

$$E_a(-y^2) = A(y)\cos\phi(y), \qquad yO_a(-y^2) = A(y)\sin\phi(y)$$

at $y = s_j$ gives

$$\frac{O_a(-s_j^2)}{E_a{}'(-s_j^2)} = \frac{2}{\phi'(s_j)} > 0.$$

Since $E_a$ and $O_a$ have the same degree, partial fractions give

$$\frac{O_a(z)}{E_a(z)} = \frac{e_{2r+1}(a)}{e_{2r}(a)} + \sum_{j=1}^{r} \frac{\gamma_j}{z + s_j^2}, \qquad \gamma_j = \frac{O_a(-s_j^2)}{E_a{}'(-s_j^2)} > 0.$$

Every summand has negative imaginary part on $\mathbb{H}$. Hence $O_a/E_a$ maps $\mathbb{H}$ into the lower half-plane, and therefore

$$E_a(z) + wO_a(z) \neq 0, \qquad z, w \in \mathbb{H}.$$

Fix $w \in \mathbb{H}$. The polynomial $E_a(z) + wO_a(z)$ has exact degree $r$, because its leading coefficient $e_{2r}(a) + we_{2r+1}(a)$ is nonzero. Applying Lemma 3.3 in the $z$-variable, while retaining $w$ as a parameter, gives

$$P_a(z) + wB_a(z) \neq 0, \qquad z \in \mathbb{H}.$$

As this holds for every $w \in \mathbb{H}$, the quotient $B_a/P_a$ maps $\mathbb{H}$ into the closed lower half-plane.

After common factors have been cancelled, every pole is real and simple. Indeed, if a reduced pole at $x_0 \in \mathbb{R}$ had order $m \geq 2$, with leading principal part $c(z - x_0)^{-m}$, then on the upper semicircle $z = x_0 + \varepsilon e^{i\theta}$ its leading imaginary part would be

$$-c\varepsilon^{-m}\sin(m\theta),$$

which assumes both signs as $0 < \theta < \pi$. This contradicts the lower-half-plane mapping property for sufficiently small $\varepsilon$. Thus every reduced pole is simple. Approaching a pole vertically shows that its residue is positive. All poles are negative because the zeros of $P_a$ are negative, and

$$\lim_{z\to\infty} \frac{B_a(z)}{P_a(z)} = \frac{e_{2r+1}(a)}{e_{2r}(a)} > 0.$$

Consequently,

$$\frac{B_a(z)}{P_a(z)} = c_\infty + \sum_{\ell=1}^{s} \frac{c_\ell}{z + \tau_\ell}, \qquad c_\infty, c_\ell, \tau_\ell > 0. \qquad (3.4)$$

For every polynomial $p$ and every $b > 0$, the beta integral gives

$$T_{b+1/2}p(z) = \frac{\Gamma(b + 1/2)}{\Gamma(b)\Gamma(1/2)} \int_0^1 t^{b-1}\,(1-t)^{-1/2} T_b p(tz)\, dt.$$

Indeed, on the monomial $z^k$ this reduces to

$$\int_0^1 t^{b+k-1}\,(1-t)^{-1/2}\, dt = \frac{\Gamma(b+k)\Gamma(1/2)}{\Gamma(b+k+1/2)}.$$

Taking $p = O_a$ and dividing by $P_a(-x)$ yields

$$\frac{R_a(-x)}{P_a(-x)} = c_r \frac{\Gamma(b + 1/2)}{\Gamma(b)\Gamma(1/2)} \int_0^1 t^{b-1}\,(1-t)^{-1/2} \frac{B_a(-tx)}{P_a(-x)}\, dt. \qquad (3.5)$$

Expanding (3.4) at $z = -x$, write

$$\frac{B_a(-x)}{P_a(-x)} = \sum_{q=0}^{\infty} \eta_q\, x^q, \qquad \eta_q \geq 0, \qquad \eta_0 = e_1(a) > 0.$$

Also write, with multiplicity,

$$P_a(z) = \prod_{j=1}^{r} \left(1 + \frac{z}{\rho_j}\right), \qquad \rho_j > 0.$$

For $0 < t < 1$,

$$\frac{B_a(-tx)}{P_a(-x)} = \frac{B_a(-tx)}{P_a(-tx)}\frac{P_a(-tx)}{P_a(-x)},$$

where

$$\frac{B_a(-tx)}{P_a(-tx)} = \sum_{q=0}^{\infty} \eta_q\, t^q x^q$$

and

$$\frac{1 - tx/\rho_j}{1 - x/\rho_j} = 1 + (1-t)\sum_{m=1}^{\infty} \rho_j^{-m}\, x^m.$$

For $\mathbf{m} = (m_1, \dots, m_r) \in \mathbb{Z}_{\geq 0}^r$, put

$$|\mathbf{m}| = \sum_{j=1}^{r} m_j\,, \qquad \ell(\mathbf{m}) = \#\{j: m_j > 0\}.$$

On $|x| < 1/2\,\min_j \rho_j$, all these series converge absolutely and uniformly for $0 \leq t \leq 1$, while the beta weight in (3.5) is integrable. Coefficient extraction may therefore be passed through the integral, giving

$$[x^N]\frac{R_a(-x)}{P_a(-x)} = c_r \frac{\Gamma(b+1/2)}{\Gamma(b)\Gamma(1/2)} \times \sum_{q+|\mathbf{m}|=N} \eta_q \left(\prod_{j:m_j>0} \rho_j^{-m_j}\right) B\left(b+q, \ell(\mathbf{m}) + \frac{1}{2}\right) > 0. \tag{3.6}$$

Every summand in (3.6) is nonnegative. For $N = 0$, strict positivity follows from $q = 0$ and $\mathbf{m} = 0$; for $N \geq 1$, it follows from $q = 0$, $m_1 = N$, and $m_j = 0$ for $j > 1$. Equivalently,

$$(-1)^N [z^N]\frac{R_a(z)}{P_a(z)} > 0 \qquad (N \geq 0). \tag{3.7}$$

By (3.7) and the second identity in (3.2),

$$[\pi^{-N}]\kappa_{2N+1}^{(2r+1)}(a)$$

is nonzero and has sign $(-1)^N$. Together with Subsections 3.1 and 3.2, this proves the degree identity (3.1) and all its sign assertions. □

**Corollary 3.4 (extremal coefficients and their common root rate).** Let $d = 2r + 1$, and write

$$P_a(z) = \prod_{j=1}^{r}\left(1 + \frac{z}{\rho_j}\right), \qquad \rho_* := \min_{1\leq j\leq r} \rho_j.$$

For every $N \geq 1$,

$$\frac{(-1)^{N-1}}{(2N)!}[\pi^{-N}]\kappa_{2N}^{(2r+1)}(a) = \frac{1}{N}\sum_{j=1}^{r} \rho_j^{-N} > 0. \tag{3.8}$$

For every $N \geq 0$,

$$\frac{(-1)^N}{(2N+1)!}[\pi^{-N}]\kappa_{2N+1}^{(2r+1)}(a) = [x^N]\frac{R_a(-x)}{P_a(-x)} > 0, \qquad (3.9)$$

and the coefficient on the right is given explicitly by the positive Stieltjes–beta sum in (3.6). The two positive sequences in (3.8) and (3.9) have the same root-exponential rate:

$$\begin{aligned} &\lim_{N\to\infty}\left(\frac{(-1)^{N-1}}{(2N)!}[\pi^{-N}]\kappa_{2N}^{(2r+1)}(a)\right)^{1/N} \\ &= \lim_{N\to\infty}\left(\frac{(-1)^N}{(2N+1)!}[\pi^{-N}]\kappa_{2N+1}^{(2r+1)}(a)\right)^{1/N} = \rho_*^{-1}. \end{aligned} \qquad (3.10)$$

**Proof.** Formula (3.8) is the logarithmic root expansion from Subsection 3.2 together with the first identity in (3.2). Formula (3.9) is the second identity in (3.2), after replacing $z$ by $-x$, and strict positivity is (3.6).

The first limit in (3.10) follows immediately from the finite positive sum in (3.8). For the second, the rational function $R_a(-x)/P_a(-x)$ has no pole in $|x| < \rho_*$, so the Cauchy–Hadamard formula gives

$$\limsup_{N\to\infty}\left([x^N]\frac{R_a(-x)}{P_a(-x)}\right)^{1/N} \leq \rho_*^{-1}.$$

Choose an index $j_*$ with $\rho_{j_*} = \rho_*$. In the positive sum (3.6), retain only the term $q = 0$, $m_{j_*} = N$, and $m_j = 0$ for $j \neq j_*$. Since $\eta_0 = e_1(a) > 0$, this gives a constant $C(a,r) > 0$, independent of $N$, such that

$$[x^N]\frac{R_a(-x)}{P_a(-x)} \geq C(a,r)\rho_*^{-N}.$$

The matching lower bound proves the second limit. □

For $d = 3$,

$$P_a(z) = 1 + \frac{2e_2(a)}{3}z, \qquad R_a(z) = \frac{e_1(a)}{2} + \frac{e_3(a)}{4}z, \qquad \rho_* = \frac{3}{2e_2(a)}.$$

Thus (3.8) becomes

$$[\pi^{-N}]\kappa_{2N}^{(3)}(a) = (2N)!\,\frac{(-1)^{N-1}}{N}\left(\frac{2e_2(a)}{3}\right)^N,$$

while, for $N \geq 1$, (3.9) becomes

$$[\pi^{-N}]\kappa_{2N+1}^{(3)}(a) = (2N+1)!\,(-1)^{N-1}\left(\frac{2e_2(a)}{3}\right)^{N-1}\frac{3e_3(a) - 4e_1(a)e_2(a)}{12}.$$

These are the one-root specializations underlying the three-dimensional degree theorem in the companion paper.

In dimension three, Theorem 3.1 recovers the cumulant-degree theorem of our companion paper [1]. Its new content is the complete parity dichotomy in every dimension, including the even-dimensional case, every odd dimension above three, and the Stieltjes positivity mechanism that prevents cancellation in every odd order.

### *3.4. Reconstruction and generic sharpness*

The parity theorem concerns the arithmetic structure of the width cumulants. The same Gaussian factorization gives a direct reconstruction theorem.

Let

$$H_i = |G_i|, \qquad S_a = \sum_{i=1}^{d} a_i H_i,$$

where the $H_i$ are independent standard half-normal variables. By Lemma 2.2,

$$\mathbb{E}S_a^n = 2^{n/2} \frac{\Gamma\big((d+n)/2\big)}{\Gamma(d/2)} m_n(a).$$

Hence the first $d$ moments of $W_a(U)$ determine the first $d$ moments, and therefore the first $d$ cumulants, of $S_a$.

Independence and homogeneity of cumulants give

$$\kappa_n(S_a) = \kappa_n(H_1) \sum_{i=1}^{d} a_i^n . \qquad (3.11)$$

To use (3.11), we must know that $\kappa_n(H_1) \neq 0$. Let

$$\mu_n = \mathbb{E}H_1^n = 2^{n/2} \frac{\Gamma\big((n+1)/2\big)}{\sqrt{\pi}}.$$

Put

$$y = \sqrt{\frac{2}{\pi}}.$$

The half-normal moments separate by parity as

$$\mu_{2k} = (2k-1)!!, \qquad \mu_{2k+1} = 2^k k!\, y.$$

The moment–cumulant recurrence therefore expresses $\kappa_n(H_1)$ as a polynomial in $y$ with rational coefficients. Its term of maximal $y$-degree is the all-singleton contribution

$$(-1)^{n-1}(n-1)!\, y^n,$$

and it cannot be cancelled by any other partition. Since $\pi$ is transcendental [19], so is $y$; the resulting nonzero rational polynomial cannot vanish at $y$. Thus every $\kappa_n(H_1)$ is nonzero.

**Theorem 3.5 (finite reconstruction).** The first $d$ moments of $W_a(U)$ determine the unordered side vector

$$\{a_1, \dots, a_d\}.$$

Equivalently, they determine the orthotope $K_a$ up to coordinate permutation.

**Proof.**
The first $d$ width moments determine $\kappa_n(S_a)$ for $1 \le n \le d$. Dividing (3.11) by the nonzero constant $\kappa_n(H_1)$ recovers the power sums

$$p_n(a) = \sum_{i=1}^{d} a_i^n, \qquad 1 \le n \le d.$$

Newton's identities then recover $e_1(a), \dots, e_d(a)$. The side lengths are the positive roots of

$$x^d - e_1(a)x^{d-1} + e_2(a)x^{d-2} - \cdots + (-1)^d e_d(a).$$

Hence their unordered multiset is determined. □

The number $d$ is generically sharp.

**Proposition 3.6 (generic necessity).** At every side vector $a$ with pairwise distinct positive coordinates, there exist arbitrarily close, noncongruent positive side vectors having the same first $d-1$ width moments. Consequently, fewer than $d$ moments do not generically determine an orthotope.

**Proof.**
The first $d-1$ width moments determine, and are determined by, the first $d-1$ power sums $p_1, \dots, p_{d-1}$. Newton's identities therefore fix

$$e_1, \dots, e_{d-1},$$

but not $e_d$.

Consider

$$q_\varepsilon(x) = x^d - e_1 x^{d-1} + \cdots + (-1)^{d-1} e_{d-1} x + (-1)^d (e_d + \varepsilon).$$

At $\varepsilon = 0$, its roots $a_1, \dots, a_d$ are positive and simple. By continuity of simple roots, equivalently by the implicit-function theorem applied at each root, for all sufficiently small real $\varepsilon$, the polynomial $q_\varepsilon$ still has $d$ distinct positive roots. Their first $d-1$ elementary symmetric functions, and hence their first $d-1$ power sums, agree with those of $a$, while their product differs. They therefore define a nearby noncongruent orthotope with the same first $d-1$ width moments. □

For $d = 3$, Theorem 3.5 recovers the uniqueness conclusion of our companion paper [1]. The two results have different strengths. The present argument works in every dimension and exposes the Gaussian power-sum mechanism. The companion paper, by contrast, provides explicit three-dimensional inversion formulas for the edge sum, surface area and volume, together with a realizability criterion for cumulant triples. Theorem 3.5 gives finite algebraic reconstruction but does not replace those explicit formulas or the realizability theorem.

### *3.5. Correlated Gaussian generators*

The Gaussian factorization extends beyond coordinate orthotopes. Let $v_1, \dots, v_M \in \mathbb{R}^d$, let $c_j > 0$, and define

$$Z(u) = \sum_{j=1}^{M} c_j \left|\langle v_j, u\rangle\right|, \qquad u \in S^{d-1}.$$

**Proposition 3.7.** If $U$ is uniform on $S^{d-1}$ and $G$ is standard Gaussian in $\mathbb{R}^d$, then for every $N \geq 0$,

$$\mathbb{E}Z(U)^N = \frac{\Gamma(d/2)}{2^{N/2}\Gamma\big((d+N)/2\big)}\,\mathbb{E}\left(\sum_{j=1}^{M} c_j \left|\langle v_j, G\rangle\right|\right)^N. \qquad (3.12)$$

**Proof.**
Write $G = RU$, with $R$ independent of $U$. Then

$$\sum_{j=1}^{M} c_j \left|\langle v_j, G\rangle\right| = R\,Z(U).$$

Taking $N$-th moments and dividing by

$$\mathbb{E}R^N = 2^{N/2}\frac{\Gamma\big((d+N)/2\big)}{\Gamma(d/2)}$$

gives (3.12). □

In the orthotopal case, the vectors $v_j$ are mutually orthogonal coordinate vectors. The Gaussian projections are then independent, and their absolute moments factor into one-dimensional half-normal moments, producing the power-sum identity (3.11). For general generators, the covariance data

$$\mathrm{Cov}\big(\langle v_j, G\rangle, \langle v_k, G\rangle\big) = \langle v_j, v_k\rangle$$

introduce mixed absolute moments governed by the full Gram matrix. The corresponding sign-cone, or hyperplane-arrangement, decomposition remains exact, but the independent half-normal factorization is lost. This Gram-matrix coupling is precisely what prevents the moment-closure argument (3.11) from extending directly to general zonotopes. Within this factorization route, orthogonal generators are exactly the case in which the Gram matrix is diagonal and Newton's identities close the reconstruction; no global non-reconstructibility claim is made for special correlated generator systems.

Proposition 3.7 will be used later when the width framework is extended from orthotopes to zonotopal and linear–quadratic models. Before that extension, Section 4 applies the orthotopal reconstruction mechanism mode by mode to the spectral family generated by a fixed box. Sections 5 and 6 then return to the discrete decomposition (2.4), so that the continuous reconstruction theory and the spectral analysis meet through a common family of geometric data rather than through subset indexing alone.

# 4. Auxiliary width laws and spectral data

Section 3 established that finitely many moments of a single orthotopal width law determine its side vector. We now place that continuous rigidity mechanism inside the separated spectral theory of an orthotope. This is also the precise point of continuity with our companion paper [1]: in dimension three, the width formulas established there apply to every nondegenerate auxiliary member of the mode-indexed family introduced below.

### *4.1. The mode-indexed family*

Let

$$\mathcal{J}_D = \mathbb{Z}^d_{>0}, \qquad \mathcal{J}_N = \mathbb{Z}^d_{\geq 0} \setminus \mathbf{0}$$

denote the Dirichlet and nonzero Neumann index sets. For every nonzero $n = (n_1, \dots, n_d) \in \mathbb{Z}^d_{\geq 0}$, define

$$b^{(n)}(a) = \left(\frac{n_1}{a_1}, \dots, \frac{n_d}{a_d}\right),$$

$$W_{n,a}(u) = \sum_{i=1}^{d} \frac{n_i}{a_i} |u_i|, \qquad u \in S^{d-1}, \tag{4.1}$$

$$\mu_{n,a} = \mathcal{L}\left(W_{n,a}(U)\right), \qquad M_{n,a} = \operatorname{maxsupp}\mu_{n,a}.$$

When $n \in \mathcal{J}_D$, $W_{n,a}$ is the width function of the full-dimensional auxiliary orthotope with side vector $b^{(n)}(a)$. If some coordinates of $n$ vanish, it is the ambient spherical width function of the corresponding embedded lower-dimensional orthotope. Equivalently, its law is the limit of full-dimensional laws as the vanishing auxiliary side lengths tend to zero. The moment formulas and Gaussian factorization of Section 2 extend to these degenerate laws by polynomial continuity; the non-cancellation conclusions of Theorem 3.1 are not asserted when positivity of every auxiliary side is lost. These degenerate members are not discarded in spectral aggregation: they are precisely the lower-dimensional coordinate-stratum contributions isolated by Lemma 2.4.

**Proposition 4.1 (spectral endpoint identity).** For every nonzero $n \in \mathbb{Z}^d_{\geq 0}$,

$$M_{n,a} = \left(\sum_{i=1}^{d} \frac{n_i^2}{a_i^2}\right)^{1/2}.$$

Moreover,

$$\lambda_n^D(a) = \pi^2 M_{n,a}^2 \quad (n \in \mathcal{J}_D), \qquad \lambda_n^N(a) = \pi^2 M_{n,a}^2 \quad (n \in \mathcal{J}_N). \tag{4.2}$$

The zero Neumann eigenvalue corresponds to $n = 0$ and is excluded from the auxiliary family.

**Proof.**
For $u \in S^{d-1}$, Cauchy–Schwarz gives

$$W_{n,a}(u) \leq \left(\sum_{i=1}^{d} \frac{n_i^2}{a_i^2}\right)^{1/2} \left(\sum_{i=1}^{d} |u_i|^2\right)^{1/2} = \left(\sum_{i=1}^{d} \frac{n_i^2}{a_i^2}\right)^{1/2}.$$

Equality is attained by taking

$$|u_i| = \frac{n_i/a_i}{\left(\sum_j n_j^2/a_j^2\right)^{1/2}}.$$

Thus the first assertion follows.

Translation does not change the spectrum, so the centred orthotope $K_a$ may be replaced by

$$Q_a = \prod_{i=1}^{d} (0, a_i).$$

Separation of variables gives the Dirichlet eigenfunctions

$$x \mapsto \prod_{i=1}^{d} \sin\left(\frac{\pi n_i x_i}{a_i}\right), \qquad n \in \mathcal{I}_D,$$

and the Neumann eigenfunctions

$$x \mapsto \prod_{i=1}^{d} \cos\left(\frac{\pi n_i x_i}{a_i}\right), \qquad n \in \mathbb{Z}_{\geq 0}^{d}.$$

Applying $-\Delta$ multiplies either product by

$$\pi^2 \sum_{i=1}^{d} \frac{n_i^2}{a_i^2}.$$

Comparison with the endpoint already obtained proves (4.2). □

Thus a single orthotope generates the infinite indexed family

$$\{\mu_{n,a}\}_{n \in \mathcal{I}_D} \qquad \text{or} \qquad \{\mu_{n,a}\}_{n \in \mathcal{I}_N},$$

and its separated eigenvalues are the squared upper endpoints of those laws, scaled by $\pi^2$.

The construction is scale covariant. For $c > 0$,

$$W_{n,ca} = c^{-1} W_{n,a}, \qquad M_{n,ca} = c^{-1} M_{n,a}, \qquad \lambda_n(ca) = c^{-2} \lambda_n(a).$$

Thus an overall dilation separates cleanly from the side ratios; the latter govern the shape of every normalized auxiliary law.

### *4.2. The information hierarchy*

Proposition 4.1 passes through four distinct levels of data:

$$\{W_{n,a}(u)\}_{n,u} \to \{\mu_{n,a}\}_n \to \{M_{n,a}\}_n \to \{\pi^2 M_{n,a}^2\}_{\text{multiset}}.$$

Each arrow is a specified forgetting operation. Passing from the direction-labelled functions to the indexed laws forgets which direction produces each value. Passing from an indexed law to its endpoint discards all interior distributional data, including its moments, cumulants and singular structure. Passing from labelled endpoints to the ordinary spectrum forgets the modal labels while retaining multiplicities.

The endpoint reduction is genuinely non-injective. On $S^1$, consider the positive coefficient vectors

$$b = \left(\frac{\sqrt{3}}{2}, \frac{1}{2}\right), \qquad c = \left(\frac{1}{\sqrt{2}}, \frac{1}{\sqrt{2}}\right).$$

Both have Euclidean norm 1, and hence both width laws have upper endpoint 1. Their means, however, are

$$\mathbb{E}\, W_b(U) = \frac{\sqrt{3}+1}{\pi}, \qquad \mathbb{E}\, W_c(U) = \frac{2\sqrt{2}}{\pi},$$

which are unequal. An endpoint therefore does not determine the law whose endpoint it is.

Accordingly, (4.2) does not identify the interior of an auxiliary law as a per-mode spectral observable. It says exactly that the eigenvalue records its squared upper endpoint. A separate inverse-spectral theorem may nevertheless recover the side vector $a$ from an ordinary unlabelled spectrum, after which the indexed family can be regenerated from (4.1). This is an indirect reconstruction through the recovered geometry, not a consequence of the endpoint map itself. Recent index-free studies treat planar rectangles and three-dimensional boxes by asymptotic spectral data [22, 23]. In dimension three, the latter route extracts the geometric coefficients that determine edge sum, surface area and volume, and then recovers the edges algebraically. The companion paper reaches the same elementary symmetric data from width cumulants. The inputs and observability questions are different, even though both reconstructions terminate at the same cubic.

### *4.3. Rigidity at three levels*

Within the orthotope class, explicit inverse statements are available before the modal labels are discarded.

**Proposition 4.2 (indexed-law and labelled-endpoint rigidity).** Let $a \in (0,\infty)^d$.

1. The first $d$ moments of the single indexed law $\mu_{\mathbf{1},a}$, where $\mathbf{1} = (1,\ldots,1)$, determine the unordered side vector $a$.
2. The $d$ labelled Neumann endpoints corresponding to $\varepsilon_1,\ldots,\varepsilon_d$, where $\varepsilon_i$ is the $i$-th standard basis vector, determine the ordered side vector.
3. The $d+1$ labelled Dirichlet eigenvalues corresponding to $\mathbf{1}$ and $\mathbf{1}+\varepsilon_i$, $1 \le i \le d$, determine the ordered side vector.

More explicitly,

$$a_i = \frac{1}{M_{\varepsilon_i,a}} = \frac{\pi}{\sqrt{\lambda^N_{\varepsilon_i}(a)}}, \qquad a_i = \pi\sqrt{\frac{3}{\lambda^D_{\mathbf{1}+\varepsilon_i}(a) - \lambda^D_{\mathbf{1}}(a)}}. \qquad (4.3)$$

**Proof.**
The law $\mu_{\mathbf{1},a}$ is the width law of the auxiliary orthotope with side vector

$$\left(\frac{1}{a_1},\ldots,\frac{1}{a_d}\right).$$

Theorem 3.5 therefore recovers the unordered reciprocal side vector from its first $d$ moments, and hence recovers the unordered vector $a$.

For the Neumann coordinate mode $\varepsilon_i$, (4.2) gives

$$M_{\varepsilon_i,a} = \frac{1}{a_i}, \qquad \lambda^N_{\varepsilon_i}(a) = \frac{\pi^2}{a_i^2}.$$

For the Dirichlet indices $\mathbf{1}$ and $\mathbf{1}+\varepsilon_i$,

$$\lambda^D_{\mathbf{1}+\varepsilon_i}(a) - \lambda^D_{\mathbf{1}}(a) = \frac{3\pi^2}{a_i^2}.$$

These identities prove (4.3). □

Proposition 4.2 connects the rigidity theorem of Section 3 directly to the spectral family: one full indexed law reconstructs the orthotope through its moments, whereas finitely many labelled endpoints reconstruct it even before the interior laws are used. Neither assertion solves the ordinary unlabelled inverse-spectral problem, because the modal identification required by (4.3) is then absent.

In dimension two, the first two ordered Dirichlet eigenvalues already determine the rectangle without modal labels. If $a_1 \geq a_2$, then the second eigenvalue is associated with $(2,1)$, so its difference from the ground eigenvalue is $3\pi^2/a_1^2$; the ground eigenvalue then determines $a_2$. Proposition 4.2 has a different role: it identifies exactly where each labelled level of the information hierarchy becomes sufficient in every dimension.

Ordinary unlabelled spectral data also give a finite reconstruction when finite data means finitely many complete-spectrum invariants rather than finitely many individual eigenvalues. The mechanism is the classical factorization of an orthotopal heat trace into one-dimensional theta series. Its leading coefficient belongs to Weyl's law [30], while its boundary coefficient belongs to the two-term spectral asymptotic tradition [14]; the point recorded here is the exact all-dimensional orthotope inversion, extending the three-dimensional use of the same geometric coefficients in [23].

**Proposition 4.3 (finite heat-invariant reconstruction).** The $d$ nonconstant algebraic coefficients in the small-time Dirichlet heat trace of an orthotope determine its unordered side vector. The same holds for Neumann conditions.

**Proof.** In one dimension, periodize

$$g(x) = e^{-\pi^2 t x^2/a^2}.$$

The Gaussian decay gives uniform convergence of the periodized series and its derivatives. If

$$I_\alpha(\xi) = \int_{\mathbb{R}} e^{-\alpha x^2}\, e^{-2\pi i \xi x}\, dx,$$

differentiation followed by integration by parts gives

$$I_\alpha{}'(\xi) = -\frac{2\pi^2 \xi}{\alpha} I_\alpha(\xi), \qquad I_\alpha(0) = \sqrt{\frac{\pi}{\alpha}}.$$

Solving this differential equation and evaluating the absolutely convergent Fourier series at zero yields

$$\sum_{n\in\mathbb{Z}} e^{-\pi^2 t n^2/a^2} = \frac{a}{\sqrt{\pi t}} \sum_{k\in\mathbb{Z}} e^{-a^2 k^2/t}.$$

Consequently,

$$\sum_{n=1}^{\infty} e^{-\pi^2 t n^2/a^2} = \frac{a}{2\sqrt{\pi t}} - \frac{1}{2} + O_a\left(t^{-1/2} e^{-a^2/t}\right).$$

Multiplication in the $d$ coordinates gives

$$Z_D(t;a) = \frac{1}{2^d}\sum_{k=0}^{d}(-1)^{d-k}\, e_k(a)(\pi t)^{-k/2} + O_a\big(t^{-d/2}e^{-c_a/t}\big). \qquad (4.4)$$

Thus the $d$ nonconstant coefficients in (4.4) recover $e_1(a), \dots, e_d(a)$. The side lengths are the positive roots of

$$z^d - e_1(a)z^{d-1} + e_2(a)z^{d-2} - \cdots + (-1)^d e_d(a). \qquad (4.5)$$

For Neumann conditions the one-dimensional constant is $+1/2$, so the same elementary symmetric functions occur with positive signs. The polynomial (4.5) then gives the same reconstruction. □

Proposition 4.3 concerns scalar invariants extracted from the complete unlabelled spectrum. It is not reconstruction from finitely many individual eigenvalues or from finitely many heat-trace values at fixed positive times.

The next section studies what remains when individual modes are aggregated into spectral populations. For an even summable lattice weight, Lemma 2.4 resolves the aggregate exactly into its full-dimensional positive-index term and its lower-dimensional zero-coordinate strata. Applied to cutoff, Riesz and heat weights constructed from (4.2), it becomes the summation mechanism connecting the auxiliary family to global spectral populations; Section 6 will retain the first non-universal stratum and extract its face-scale geometric memory.

# 5. Spectral universality, conditional optimization, and all-scale balance

Section 4 distinguished an indexed auxiliary law from its endpoint and from the ordinary unlabelled spectrum. We now retain the separated Dirichlet labels and ask two different questions. First, which fixed-volume orthotope minimizes a prescribed labelled spectral population? Second, what distribution of directional partitions is produced when all labels below a growing spectral cutoff are aggregated? The first problem is shape-selective but population-dependent; the second has a universal leading limit.

## *5.1. Labelled partitions and conditional optimization*

For $n \in \mathcal{I}_D$, define

$$q_i(n;a) = \frac{n_i^2/a_i^2}{\sum_{j=1}^{d} n_j^2/a_j^2}, \qquad q(n;a) = \big(q_1(n;a), \dots, q_d(n;a)\big). \qquad (5.1)$$

Thus $q(n;a)$ belongs to the open simplex and records the coordinate contributions to $M_{n,a}^2$. It is attached to the canonical separated label $n$; at a multiple eigenvalue it is not invariant under an arbitrary change of eigenbasis.

Equivalently, the positive-orthant maximizing direction in Proposition 4.1 is

$$u_{n,a}^* = \frac{b^{(n)}(a)}{M_{n,a}} = \frac{1}{M_{n,a}}\Big(\frac{n_1}{a_1}, \dots, \frac{n_d}{a_d}\Big),$$

and

$$q_i(n;a) = \big(u_{n,a,i}^*\big)^2.$$

Thus the partition is the squared-coordinate vector of the direction that maximizes the corresponding auxiliary width.

Fix $V > 0$ and write

$$\mathcal{C}_V = \{a \in (0,\infty)^d : \prod_{i=1}^{d} a_i = V\}.$$

A prescribed labelled population is a probability family $\omega = (\omega_n)_{n\in\mathcal{I}_D}$ whose support and weights are transported unchanged as $a$ varies, and such that

$$Q_i = \sum_{n\in\mathcal{I}_D} \omega_n\, n_i^2 < \infty \qquad (1 \le i \le d).$$

Since $n_i \ge 1$ on $\mathcal{I}_D$, one has $Q_i \ge 1$. The population-averaged eigenvalue is

$$\mathcal{E}_\omega(a) = \sum_n \omega_n\, \lambda_n^D(a) = \pi^2 \sum_{i=1}^{d} \frac{Q_i}{a_i^2}.$$

This shape-independence is a structural hypothesis, not a technical qualification. It excludes cutoff and heat populations, whose membership or weights change with the spectrum, and it excludes the reordered $k$-th eigenvalue problem, whose separated label may change with the shape. In what follows, an ordered side vector is coordinate-labelled rather than identified modulo coordinate permutations.

**Theorem 5.1 (labelled non-dominance and conditional realization).** The following statements hold on $\mathcal{C}_V$.

1. For every fixed label $n \in \mathcal{I}_D$,

$$\lambda_n^D(a) \ge d\pi^2 \left(\frac{\prod_i n_i}{V}\right)^{2/d},$$

with equality at the unique ordered side vector

$$a_i = V^{1/d} \frac{n_i}{\left(\prod_j n_j\right)^{1/d}}.$$

2. If $a, b \in \mathcal{C}_V$ are distinct ordered side vectors, then there are labels $n, m \in \mathcal{I}_D$ such that

$$\lambda_n^D(a) < \lambda_n^D(b), \qquad \lambda_m^D(a) > \lambda_m^D(b).$$

3. Every prescribed labelled population has the unique minimizer

$$a_i^*(\omega) = V^{1/d} \frac{\sqrt{Q_i}}{\left(\prod_j Q_j\right)^{1/(2d)}}. \qquad (5.2)$$

4. Conversely, every $a \in \mathcal{C}_V$ is the unique minimizer for a population supported on at most $d$ labels.

**Proof.**
For a fixed label, put $x_i = n_i^2/a_i^2$. Since

$$\prod_{i=1}^{d} x_i = \frac{(\prod_i n_i)^2}{V^2},$$

the arithmetic–geometric mean inequality gives the first assertion. Equality holds exactly when all $x_i$ are equal, and the volume constraint then gives the displayed side vector.

For the second assertion, equality of the products and $a \neq b$ imply that $a_i > b_i$ for some $i$ and $a_j < b_j$ for some $j$. Let $n^{(i)}(N)$ have entry $N$ in coordinate $i$ and entry 1 elsewhere. Then

$$\frac{\lambda^D_{n^{(i)}(N)}(a) - \lambda^D_{n^{(i)}(N)}(b)}{\pi^2} = N^2\left(a_i^{-2} - b_i^{-2}\right) + \sum_{k \neq i} \left(a_k^{-2} - b_k^{-2}\right).$$

The leading coefficient is negative, so the difference is negative for all sufficiently large $N$. The coordinate $j$ gives the reverse inequality.

For a population, apply the same arithmetic–geometric mean inequality to $Q_i/a_i^2$. Equality requires $a_i \propto \sqrt{Q_i}$, and the volume constraint gives (5.2).

It remains to prove realization. Given $a$, choose an integer $M > 1$ such that

$$(M^2 + d - 1)\min_i a_i^2 > \sum_{j=1}^{d} a_j^2.$$

For $1 \leq i \leq d$, let $n^{(i)}$ have entry $M$ in coordinate $i$ and entry 1 elsewhere. Put

$$C = \frac{M^2 + d - 1}{\sum_j a_j^2}, \qquad \omega_i = \frac{C a_i^2 - 1}{M^2 - 1}.$$

The choice of $M$ gives $\omega_i > 0$, while direct summation gives $\sum_i \omega_i = 1$. The second label moments of this population are

$$Q_i = 1 + (M^2 - 1)\omega_i = C a_i^2.$$

Formula (5.2) therefore returns the prescribed side vector $a$, and uniqueness follows from the equality case of the arithmetic–geometric mean inequality. □

Theorem 5.1 concerns prescribed shape-independent populations on separated labels. In the ordered-eigenvalue problem, the labels are resorted as the orthotope varies, so no fixed population represents the $k$-th eigenvalue. For rectangles, asymptotic optimization results select the square in the Dirichlet and Neumann problems [3, 26]; the corresponding three-dimensional Dirichlet result selects the cube [27]. Gittins and Larson proved the all-dimensional Dirichlet and Neumann cuboid limits and related Riesz-mean results [10]. Those theorems concern reordered spectra and the limit $k \to \infty$. By contrast, Theorem 5.1 proves that, within a fixed-volume class, no orthotope minimizes every prescribed shape-independent labelled population, while every orthotope is uniquely optimal for some such population supported on at most $d$ labels.

### *5.2. The universal limiting partition law*

For $a \in (0,\infty)^d$, set

$$V(a) = \prod_{i=1}^{d} a_i.$$

For $R > 0$, let

$$\mathcal{N}_R(a) = \{n \in \mathbb{Z}_{>0}^d : \sum_{i=1}^{d} \frac{n_i^2}{a_i^2} \le R^2\}, \qquad N_R(a) = |\mathcal{N}_R(a)|,$$

and, whenever $N_R(a) > 0$, form the empirical probability measures

$$\nu_{R,a} = \frac{1}{N_R(a)} \sum_{n \in \mathcal{N}_R(a)} \delta_{q(n;a)}, \qquad \eta_{R,a} = \frac{1}{N_R(a)} \sum_{n \in \mathcal{N}_R(a)} \delta_{u_{n,a}^*}.$$

Directional lattice-point asymptotics in expanding ellipsoidal and anisotropic domains belong to an established theory; see, for example, [18]. The result below does not claim the underlying directional equidistribution as new. Its purpose is to identify the squared-coordinate pushforward with the same Dirichlet law as in Proposition 2.1, to isolate the coordinate strata through Lemma 2.4, and to derive the geometric consequences used below and in Section 6.

**Theorem 5.2 (spectral partition and maximizing-direction universality).** For every $a \in (0,\infty)^d$, the partitions defined by (5.1) have empirical measures $\nu_{R,a}$ that converge weakly, as $R \to \infty$, to the Dirichlet law with parameter vector $(1/2, \dots, 1/2)$. Equivalently, for every continuous $\Phi$ on $\overline{\Delta}_{d-1}$,

$$\lim_{R\to\infty} \frac{1}{N_R(a)} \sum_{n \in \mathcal{N}_R(a)} \Phi\big(q(n;a)\big) = \frac{\Gamma(d/2)}{\Gamma(1/2)^d} \int_{\Delta_{d-1}} \Phi(x) \prod_{i=1}^{d} x_i^{-1/2}\, dx_1 \cdots dx_{d-1}. \tag{5.3}$$

The limit is independent of the side vector. Equivalently, $\eta_{R,a}$ converges weakly to normalized surface measure on the positive spherical orthant $S_+^{d-1}$.

**Proof.**
Extend the cutoff weight coordinatewise evenly to $\mathbb{Z}^d$, assigning any bounded value at the origin. Lemma 2.4 resolves its positive-index sum into full-lattice sums over the coordinate subspaces. On a stratum with free-coordinate set $J$, the cutoff contains at most

$$\prod_{j \in J} \big(2a_j R + 1\big) = O_a\big(R^{|J|}\big)$$

lattice points. Hence every proper stratum is $O_a\big(R^{d-1}\big)$. The zero-coordinate members are retained exactly by (2.4); their disappearance from the normalized leading limit will follow from the positive full-dimensional point-count asymptotic below.

More explicitly, if $f_R$ denotes the even cutoff weight with test function $\Phi$, then

$$\sum_{n\in\mathcal{N}_R(a)} \Phi\left(q(n;a)\right) = 2^{-d}L_{[d]}(f_R) + O_{a,\Phi}\left(R^{d-1}\right),$$

and the same identity with $\Phi \equiv 1$ holds for $N_R(a)$.

In the full-dimensional term, set

$$z_i = \frac{n_i}{a_i R}.$$

The lattice mesh has cell volume $\left(V(a)R^d\right)^{-1}$. The function

$$z \mapsto \mathbf{1}_{\{|z|\le 1\}}\Phi\left(\frac{z_1^2}{|z|^2}, \dots, \frac{z_d^2}{|z|^2}\right)$$

may be defined arbitrarily at the origin; it is bounded and Riemann integrable because its discontinuity set has Lebesgue measure zero. The lattice sums therefore converge to the corresponding integral over the unit ball. In particular,

$$N_R(a) = \frac{V(a)\kappa_d}{2^d}R^d + o_a\left(R^d\right),$$

where $\kappa_d$ is the volume of the unit ball. Every proper coordinate stratum therefore disappears after normalization by $N_R(a)$. Taking the ratio with the corresponding test-function asymptotic, radial integration cancels and leaves normalized surface measure on $S^{d-1}$. Proposition 2.1 identifies its squared-coordinate pushforward with the measure on the right-hand side of (5.3). Since the coordinatewise square map is a homeomorphism from $S_+^{d-1}$ onto $\overline{\Delta}_{d-1}$, the assertion for $\eta_{R,a}$ follows as well. □

At leading order, Lemma 2.4 isolates every proper coordinate stratum and the normalization then suppresses it. Section 6 will instead retain the strata with $d-1$ free coordinates and compute their first non-universal contribution.

**Corollary 5.3 (modal-width correspondence).** For $n \in \mathcal{I}_D$, put

$$Y_{n,a} = W_a\left(u_{n,a}^*\right) = \sum_{i=1}^{d} a_i\, u_{n,a,i}^* = \frac{\sum_{i=1}^{d} n_i}{M_{n,a}}.$$

Then the empirical laws of $Y_{n,a}$ over $\mathcal{N}_R(a)$ converge weakly to the width law of $K_a$:

$$\frac{1}{N_R(a)}\sum_{n\in\mathcal{N}_R(a)} \delta_{Y_{n,a}} \;\Rightarrow\; \mathcal{L}\left(W_a(U)\right) \qquad (R\to\infty).$$

**Proof.**
By Theorem 5.2, the partitions $q(n;a)$ converge to the law of $\left(U_1^2, \dots, U_d^2\right)$. Apply the continuous map

$$x \mapsto \sum_{i=1}^{d} a_i\sqrt{x_i}$$

on $\overline{\Delta}_{d-1}$. Since $u_{n,a,i}^* = \sqrt{q_i(n;a)}$, its value at $q(n;a)$ is $Y_{n,a}$; its value at $\left(U_1^2, \dots, U_d^2\right)$ is $\sum_i a_i\,|U_i| = W_a(U)$. □

Corollary 5.3 is a correspondence for separated labels and their maximizing directions. It does not assert that the ordinary unlabelled spectrum alone contains these direction data.

For the quadratic balance defect

$$\mathcal{B}_d(x) = d\sum_{i=1}^{d} x_i^2 - 1,$$

Proposition 2.1 gives

$$\mathbb{E}X_i^2 = \frac{(1/2)(3/2)}{(d/2)(d/2+1)} = \frac{3}{d(d+2)}.$$

Consequently,

$$\lim_{R\to\infty}\frac{1}{N_R(a)}\sum_{n\in\mathcal{N}_R(a)}\mathcal{B}_d\left(q(n;a)\right) = \frac{2(d-1)}{d+2}.$$

The limiting individual partitions therefore do not concentrate at $(1/d,\dots,1/d)$. Nevertheless, if

$$A_i(R;a) = \pi^2\sum_{n\in\mathcal{N}_R(a)}\frac{n_i^2}{a_i^2},$$

the same Riemann-sum argument and coordinate symmetry give

$$A_i(R;a) = \frac{\pi^2 V(a)\kappa_d}{2^d(d+2)}R^{d+2} + o_a\left(R^{d+2}\right),$$

where $\kappa_d$ is the volume of the unit ball. Thus

$$\frac{A_i(R;a)}{\sum_j A_j\,(R;a)} \to \frac{1}{d}.$$

Individual modal imbalance persists under a universal law, while the aggregate balances by spherical symmetry. Hence any continuous functional depending only on the limiting partition law takes the same value for every orthotope and cannot select a strict optimal ratio. This conclusion concerns the leading limiting law, not finite cutoffs.

### *5.3. Coordinate-resolved rigidity and exact spectral-scale balance*

The cube is nevertheless distinguished by exact balance at every spectral scale. For $t > 0$, define the Dirichlet heat population and its coordinate mean contributions by

$$Z_a(t) = \sum_{n\in\mathcal{I}_D} e^{-t\lambda_n^D(a)}, \qquad U_i(t;a) = \frac{1}{Z_a(t)}\sum_{n\in\mathcal{I}_D}\frac{\pi^2 n_i^2}{a_i^2}e^{-t\lambda_n^D(a)}.$$

**Theorem 5.4 (coordinate-resolved rigidity and all-scale balance).** For a Dirichlet orthotope, define its coordinate-resolved ground-state contribution vector by

$$c(a) = \left(\frac{\pi^2}{a_1^2},\dots,\frac{\pi^2}{a_d^2}\right).$$

This vector determines the ordered side vector through

$$a_i = \frac{\pi}{\sqrt{c_i(a)}}.$$

Moreover, the following conditions are equivalent:

1. $a_1 = \cdots = a_d$;
2. $c_1(a) = \cdots = c_d(a)$;
3. for every $R > 0$ with $N_R(a) > 0$, the cutoff sums $A_1(R; a), \dots, A_d(R; a)$ are equal;
4. for every $t > 0$, the canonical mean coordinate contributions $U_1(t; a), \dots, U_d(t; a)$ are equal.

Every orthotope satisfies

$$\frac{A_i(R; a)}{\sum_j A_j\,(R; a)} \to \frac{1}{d} \quad (R \to \infty), \qquad \frac{U_i(t; a)}{\sum_j U_j\,(t; a)} \to \frac{1}{d} \quad (t \downarrow 0).$$

**Proof.**
The recovery formula for $a_i$ is immediate from the definition of $c_i(a)$. Conditions 1 and 2 are therefore equivalent.

If $a_1 = \cdots = a_d$, the label set and the cutoff and heat weights appearing in the theorem are invariant under coordinate permutations, so the equalities in conditions 3 and 4 hold at every scale.

Conversely, choose $R$ immediately above the ground-state radius and below the next distinct radius. Only the label **1** contributes, and

$$A_i(R; a) = \frac{\pi^2}{a_i^2} = c_i(a).$$

Condition 3 therefore implies condition 2.

For the canonical population, separation of variables gives

$$Z_a(t) = \sum_{n \in \mathcal{I}_D} e^{-t\lambda_n^D(a)} = \prod_{i=1}^{d} \vartheta_i\,(t), \qquad \vartheta_i(t) = \sum_{m=1}^{\infty} e^{-t\pi^2 m^2 / a_i^2},$$

and

$$U_i(t; a) = \frac{\pi^2}{a_i^2} \frac{\sum_{m \ge 1} m^2\, e^{-t\pi^2 m^2 / a_i^2}}{\sum_{m \ge 1} e^{-t\pi^2 m^2 / a_i^2}}. \qquad (5.4)$$

As $t \to \infty$, the term $m = 1$ dominates both series in (5.4), so

$$U_i(t; a) \to \frac{\pi^2}{a_i^2}.$$

Condition 4 therefore implies condition 2.

Finally, with $s = t\pi^2 / a_i^2$, elementary one-dimensional Riemann sums give

$$\sqrt{s}\sum_{m\geq 1} e^{-sm^2} \to \int_0^\infty e^{-x^2}\,dx,$$

and

$$s^{3/2}\sum_{m\geq 1} m^2\,e^{-sm^2} \to \int_0^\infty x^2\,e^{-x^2}\,dx = \frac{1}{2}\int_0^\infty e^{-x^2}\,dx.$$

Substitution in (5.4) yields $U_i(t;a) \sim (2t)^{-1}$, independently of $i$ and $a$. The cutoff asymptotic was proved in Subsection 5.2. □

The proof contains two sharper finite-scale statements. If $\rho_1(a) = M_{\mathbf{1},a}$ and $\rho_2(a) = \min_{n\neq\mathbf{1}} M_{n,a}$, then every single cutoff $\rho_1(a) \leq R < \rho_2(a)$ gives

$$\big(A_1(R;a),\ldots,A_d(R;a)\big) = \left(\frac{\pi^2}{a_1^2},\ldots,\frac{\pi^2}{a_d^2}\right).$$

Thus one cutoff vector in the ground-state window recovers the ordered side vector. Likewise, the vector limit

$$\lim_{t\to\infty}\big(U_1(t;a),\ldots,U_d(t;a)\big) = \left(\frac{\pi^2}{a_1^2},\ldots,\frac{\pi^2}{a_d^2}\right)$$

recovers it, including along any sequence $t_k \to \infty$. All-scale equality characterizes the cube, but all-scale data are not required for reconstruction.

The vector $c(a)$ is coordinate-resolved enriched data: the scalar ground eigenvalue is only $\sum_i c_i\,(a)$ and does not by itself determine the side vector. Theorem 5.2 identifies the universal leading partition law, Corollary 5.3 recovers the companion width law as the high-energy empirical law of the modal maximizing directions, and Theorem 5.4 identifies the cube as the unique orthotope exhibiting exact balance at every cutoff and every heat scale. Theorem 5.1 shows that, within a fixed-volume class, no orthotope minimizes every prescribed shape-independent labelled population, while every orthotope is uniquely optimal for some finite such population. These statements concern different data levels and are mutually compatible. Section 6 will retain the first coordinate-stratum correction suppressed by the leading limit and prove that it carries recoverable face-scale shape information.

# 6. Sharp spectral asymptotics, boundary tomography, and geometric memory

The universal law in Theorem 5.2 comes from the full-dimensional lattice stratum. Geometry first re-enters through the coordinate faces. The purpose of this section is to retain that contribution in the raw, unsmoothed cutoff, to identify its complete simplex-valued form, and then to determine how far the expansion can be continued.

The outcome has two parts. First, the normalized raw directional measure has a universal bulk law and an explicit bulk–face correction whose singular part recovers every side length. Second, the uncancelled energy cutoff has jumps of edge size, so it cannot possess a universal third coefficient with a remainder smaller than edge scale. Exact Boolean differences of separable Dirichlet and Neumann populations remove this obstruction by isolating individual coordinate strata. The resulting codimension-two coefficient is genuine and again reconstructs the orthotope.

### *6.1. Raw directional measures and homogeneous discrepancy*

For $n \neq 0$, recall

$$M_{n,a} = \left( \sum_{i=1}^{d} \frac{n_i^2}{a_i^2} \right)^{1/2}, \qquad q_i(n;a) = \frac{n_i^2/a_i^2}{M_{n,a}^2}.$$

Thus $q(n;a) \in \Delta_{d-1}$ is the squared-coordinate vector of the maximizing direction of the auxiliary width law. For $s \in \{0,2\}$, put

$$c_0 = 1, \qquad c_2 = \pi^2,$$

and define

$$\mathfrak{P}_{R,a,s}^{D} = c_s \sum_{\substack{n \in \mathbb{Z}_{>0}^d \\ M_{n,a} \le R}} M_{n,a}^s \, \delta_{q(n;a)}.$$

The Neumann measure $\mathfrak{P}_{R,a,s}^{N}$ is defined by summing over $\mathbb{Z}_{\ge 0}^d \setminus \mathbf{0}$. The case $s = 0$ is the raw counting population; $s = 2$ is the energy-weighted population.

The exponent below is the weighted fixed-ellipsoid instance of the classical positive-curvature discrepancy scale originating with Hlawka [13]. We include the direct Fourier–Poisson proof required for the present angular weights.

**Lemma 6.1 (smooth-angular homogeneous discrepancy).** Let $m \ge 2$, $\alpha \in (0,\infty)^m$, $s \in \{0,2\}$, and $\Phi \in C^\infty(\Delta_{m-1})$, smooth up to the boundary. Put

$$z(k) = \left( \frac{k_1}{\alpha_1}, \dots, \frac{k_m}{\alpha_m} \right), \qquad q(z) = \frac{(z_1^2, \dots, z_m^2)}{|z|^2}.$$

Then

$$\begin{aligned} &\sum_{k \in \mathbb{Z}^m \setminus \mathbf{0}} |z(k)|^s \, \Phi\Big(q\big(z(k)\big)\Big) \, \mathbf{1}_{\{|z(k)| \le R\}} \\ &= V(\alpha) \frac{|S^{m-1}|}{m+s} R^{m+s} \int_{\Delta_{m-1}} \Phi \, d\mu_m + O_{\alpha,\Phi}\big(R^{m+s-2+2/(m+1)}\big), \end{aligned} \tag{6.1}$$

where $V(\alpha) = \prod_i \alpha_i$ and $\mu_m = \mathrm{Dirichlet}(1/2, \dots, 1/2)$. For $m = 1$, with $\Delta_0 = \{1\}$ and $\mu_1 = \delta_1$, the same leading term holds with remainder $O_{\alpha,\Phi}(R^s)$.

**Proof.** After the diagonal change of variables, the lattice has density $V(\alpha)$. Set

$$w(x) = |x|^s \Phi\left( \frac{x_1^2}{|x|^2}, \dots, \frac{x_m^2}{|x|^2} \right).$$

For $s = 0$, first excise $|x| \leq 2\delta$. Once $R$ is large, this ball contains no nonzero point of the fixed lattice $z(\mathbb{Z}^m)$, while its zero-mode contribution is $O_\Phi(\delta^m)$. The origin itself is excluded from the sum. For bounded $R$, the same contribution is absorbed by enlarging the implicit constant. The remaining weight has the uniform oscillation bound used below. For $s = 2$, no excision is needed because the weight is continuous at the origin.

Choose a nonnegative radial even mollifier $\rho_\delta$, supported in a ball of radius $\delta$. After adding and later subtracting a constant, we may suppose that the angular factor is nonnegative. Inner and outer radial convolutions, together with the oscillation bound $O_\Phi(R^{s-1}\delta)$, give pointwise lower and upper bounds for $w\mathbf{1}_{B_R}$. Their zero modes differ by $O_\Phi(R^{m+s-1}\delta)$.

For $F = w\mathbf{1}_{B_1}$, boundary charts followed by one-dimensional integration by parts at the two stationary endpoints give

$$\left|\hat{F}(\xi)\right| \leq C_\Phi(1 + |\xi|)^{-(m+1)/2}.$$

The conic point at the origin causes no larger term. Decompose a neighbourhood of the origin into dyadic annuli. The annuli of radius at most $|\xi|^{-1}$ contribute $O(|\xi|^{-m-s})$ by absolute value. On each remaining rescaled annulus, repeated integration by parts gives the same bound. Hence the boundary decay dominates globally.

Poisson summation for the smoothed bounds gives, for the nonzero dual frequencies,

$$O_{\alpha,\Phi}\left(R^{s+(m-1)/2}\delta^{-(m-1)/2}\right).$$

Taking $\delta = R^{-(m-1)/(m+1)}$ balances this term with the boundary-layer error and gives the remainder in (6.1). Finally, polar coordinates and Proposition 2.1 yield

$$\int_{B_R} w\,(x)\,dx = \frac{|S^{m-1}|}{m+s} R^{m+s} \int_{\Delta_{m-1}} \Phi\ d\mu_m.$$

This proves (6.1). The stated $m = 1$ formula follows directly by pairing the positive and negative integers. □

The estimate is proved for the unsmoothed cutoff. Radial mollification occurs only within the upper-and-lower-bound argument and is removed before the conclusion.

### *6.2. Measure-valued universality and complete face memory*

Let $\mu_{d-1}^{(j)}$ denote the Dirichlet$(1/2, \dots, 1/2)$ law on $q_j = 0$, embedded in $\Delta_{d-1}$.

**Theorem 6.2 (raw universality and face-memory measure).** Against every test function $\Phi \in C^\infty(\Delta_{d-1})$, smooth up to the boundary,

$$\begin{aligned}\mathfrak{P}_{R,a,s}^{D} &= \frac{c_s V(a)\left|S^{d-1}\right|}{2^d(d+s)} R^{d+s}\mu_d \\ &- \frac{c_s V(a)\left|S^{d-2}\right|}{2^d(d+s-1)} R^{d+s-1} \sum_{j=1}^{d} \frac{1}{a_j}\mu_{d-1}^{(j)} \qquad (6.2) \\ &+ O_{a,\Phi}\left(R^{d+s-2+2/(d+1)}\right),\end{aligned}$$

whereas the Neumann formula has a plus sign before the face sum.

**Proof.** For a coordinatewise-even summand vanishing at the origin, Lemma 2.4 gives

$$\mathfrak{P}^D = 2^{-d} \sum_{J\subseteq[d]} (-1)^{d-|J|} L_J, \qquad \mathfrak{P}^N = 2^{-d} \sum_{J\subseteq[d]} L_J.$$

Apply Lemma 6.1 in dimension $d$ to the full stratum and in dimension $d-1$ to every coordinate face. The full stratum produces $\mu_d$. The face omitting coordinate $j$ has lattice density $V(a)/a_j$, angular law $\mu_{d-1}^{(j)}$, negative Dirichlet sign, and positive Neumann sign. Every lower stratum is $O_a\big(R^{d+s-2}\big)$. The face remainders are $O_a\big(R^{d+s-3+2/d}\big)$ for $d \geq 3$, while direct summation treats $d = 2$. These terms are absorbed by the full-dimensional remainder. □

After normalization by total mass, in the weak sense against the same class of smooth test functions,

$$\overline{\mathfrak{P}}_{R,a,s}^{D} = \mu_d + \frac{1}{R}\mathfrak{M}_{a,s} + O_{a,\Phi}\big(R^{-2+2/(d+1)}\big),$$

$$\overline{\mathfrak{P}}_{R,a,s}^{N} = \mu_d - \frac{1}{R}\mathfrak{M}_{a,s} + O_{a,\Phi}\big(R^{-2+2/(d+1)}\big), \qquad (6.3)$$

where

$$\mathfrak{M}_{a,s} = \gamma_{d,s} \sum_{j=1}^{d} \frac{1}{a_j}\Big(\mu_d - \mu_{d-1}^{(j)}\Big), \qquad \gamma_{d,s} = \frac{(d+s)\big|S^{d-2}\big|}{(d+s-1)\big|S^{d-1}\big|}. \qquad (6.4)$$

If $F_j^{\circ}$ is the relative interior of the $j$-th simplex face, then the definition (6.4) gives

$$\mathfrak{M}_{a,s}\big(F_j^{\circ}\big) = -\frac{\gamma_{d,s}}{a_j}, \qquad a_j = -\frac{\gamma_{d,s}}{\mathfrak{M}_{a,s}\big(F_j^{\circ}\big)}. \qquad (6.5)$$

The measure $\mathfrak{M}_{a,s}$ is therefore a signed bulk–face correction measure; its singular part is supported on the simplex faces. Equation (6.3) first determines this finite signed measure uniquely through its action on smooth tests. Formula (6.4) then permits evaluation on Borel sets. In particular, the relative face interiors are continuity sets for $\mathfrak{M}_{a,s}$, and (6.5) recovers the complete labelled side vector. This is an iterated limiting statement, not the evaluation of a finite-cutoff Dirichlet empirical measure on a face: the latter assigns every simplex face mass zero.

The Dirichlet and Neumann populations form an exact separator. Their half-sum cancels face memory and retains the universal bulk in (6.3); their difference cancels the bulk and doubles the face term. Taking the first simplex moment gives the corresponding tensor theorem.

For a bounded Lipschitz domain $\Omega$, put

$$\mathsf{F}(\Omega) = \int_{\partial\Omega} (I - \nu \otimes \nu)\, dS.$$

Also put

$$\mathsf{N}(\Omega) = \int_{\partial\Omega} \nu \otimes \nu\, dS.$$

This surface-normal moment is the classical rank-two interfacial Minkowski tensor; in the normalization of [17], $\mathsf{N}(\Omega) = W_1^{0,2}(\Omega)$. Consequently,

$$\mathsf{F}(\Omega) = |\partial\Omega| I - \mathsf{N}(\Omega).$$

The name *face-memory tensor* refers here to its spectral role, not to the invention of the tensor. Set

$$C_d = \frac{|S^{d-1}|}{2^d d(d+2)\pi^d}, \qquad D_d = \frac{|S^{d-2}|}{2^d (d-1)(d+1)\pi^{d-1}}.$$

Let $\{\varphi_k^\eta\}_{k\geq 1}$ be any orthonormal eigenbasis of the Dirichlet or Neumann Laplacian, with eigenvalues $\lambda_k^\eta$ repeated according to multiplicity, and define

$$\mathsf{E}_a^\eta(\Lambda) = \sum_{\lambda_k^\eta(a)\leq\Lambda} \int_{K_a} \nabla\varphi_k^\eta(x) \otimes \nabla\varphi_k^\eta(x)\, dx.$$

The sum on each eigenspace is basis-independent because it is the trace of the corresponding finite-rank spectral projector against the tensor-valued quadratic form $(f, g) \mapsto \int \nabla f \otimes \nabla g$. In the separated product basis, one-dimensional sine–cosine orthogonality makes the mixed entries vanish, and the $i$-th diagonal entry contributed by a normalized mode is $\pi^2 n_i^2 / a_i^2$. For Neumann modes, $n_i = 0$ contributes zero to that entry.

**Corollary 6.3 (projector-gradient Weyl tensor).** For Dirichlet or Neumann conditions,

$$\mathsf{E}_a^{D/N}(\Lambda) = C_d V(a) \Lambda^{d/2+1} I \mp \frac{D_d}{2} \mathsf{F}(K_a) \Lambda^{(d+1)/2} + O_a\big(\Lambda^{d/2+1/(d+1)}\big). \qquad (6.6)$$

**Proof.** Choose the separated product eigenbasis and take $s = 2$, $\Phi(q) = q_i$ in (6.2), with $R = \sqrt{\Lambda}/\pi$. Proposition 2.1 gives

$$\int q_i \, d\mu_d = \frac{1}{d}, \qquad \int q_i \, d\mu_{d-1}^{(j)} = \frac{1-\delta_{ij}}{d-1}.$$

The bulk term is therefore $C_d V(a) \Lambda^{d/2+1}$. The face term in the $i$-th diagonal entry is

$$\mp D_d V(a) \sum_{j\neq i} \frac{1}{a_j} \Lambda^{(d+1)/2} = \mp \frac{D_d}{2} \mathsf{F}(K_a)_{ii} \Lambda^{(d+1)/2}.$$

The mixed entries vanish by the product-basis orthogonality proved above, and the remainder in (6.2) becomes the remainder in (6.6). Basis independence then proves the formula for every orthonormal eigenbasis. The reconstruction below is an algebraic consequence of the displayed diagonal entries. □

Consequently,

$$\Lambda^{-(d+1)/2}\big(\mathsf{E}_a^N(\Lambda) - \mathsf{E}_a^D(\Lambda)\big) \to D_d \mathsf{F}(K_a).$$

In the principal-axis frame,

$$\mathsf{F}(K_a)_{ii} = 2V(a) \sum_{j\neq i} \frac{1}{a_j}, \qquad \mathsf{F}(K_a)_{ij} = 0 \quad (i \neq j).$$

If $f_1, \dots, f_d$ are its eigenvalues and

$$b_i = \frac{1}{2}\left(\frac{\sum_{k=1}^{d} f_k}{d-1} - f_i\right),$$

then $b_i = V(a)/a_i$, whence

$$V(a) = \left(\prod_{i=1}^{d} b_i\right)^{1/(d-1)}, \qquad a_i = \frac{V(a)}{b_i}. \qquad (6.7)$$

Moreover,

$$\mathrm{tr}\mathsf{E}_a^{\eta}(\Lambda) = \sum_{\lambda_n^{\eta}(a)\le\Lambda} \lambda_n^{\eta}(a), \qquad \mathrm{tr}\mathsf{F}(K_a) = (d-1)|\partial K_a|.$$

Thus the trace of (6.6) recovers Weyl's leading scalar coefficient and the classical second boundary coefficient, within the scalar Riesz-mean context [8, 14, 30]; the directional deviatoric tensor and its explicit inversion are the additional information retained here. Proposition 4.3 already reconstructs an orthotope from the complete heat invariants. The role of (6.7) is sharper: one tensorial boundary coefficient itself carries the complete side geometry.

Equations (6.5) and (6.7) are two projections of the same face-memory object: the former reads its simplex-face masses, while the latter reads its first tensor moment.

### *6.3. The raw edge-scale obstruction*

The face coefficient is visible because the remainder in (6.2) is smaller than face scale. The corresponding uncancelled energy cutoff cannot be continued universally through edge scale.

Put

$$N_d^+(R) = \#\{n \in \mathbb{Z}_{>0}^d : |n| \le R\}, \qquad r_d^+(m) = \#\{n \in \mathbb{Z}_{>0}^d : |n|^2 = m\}.$$

Unit-cube comparison with the positive-orthant ball gives

$$N_d^+(R) = 2^{-d}\kappa_d R^d + O_d\left(R^{d-1}\right).$$

It follows that

$$\sum_{X<m\le 2X} r_d^+(m) = 2^{-d}\kappa_d\left(2^{d/2}-1\right)X^{d/2} + O_d\left(X^{(d-1)/2}\right).$$

Hence some $m \in (X, 2X]$ satisfies $r_d^+(m) \ge c_d m^{d/2-1}$ for every sufficiently large dyadic $X$. For

$$E_d^+(R) = \sum_{\substack{n\in\mathbb{Z}_{>0}^d \\ |n|\le R}} |n|^2,$$

the jump at $R_m = \sqrt{m}$ obeys

$$\Delta E_d^+(R_m) = m r_d^+(m) \ge c_d R_m^d \qquad (6.8)$$

along an unbounded sequence.

**Theorem 6.4 (edge-scale jump obstruction).** Already for the cube, the uncancelled raw energy cutoff does not admit a universal deterministic expansion through order $R^d$ with an $o\left(R^d\right)$ remainder.

**Proof.** Suppose that

$$E_d^+(R) = C_{d+2}R^{d+2} + C_{d+1}R^{d+1} + C_dR^d + o\left(R^d\right).$$

Choose $\varepsilon_m = o(R_m^{-1})$ so that no other squared lattice radius lies between $R_m - \varepsilon_m$ and $R_m + \varepsilon_m$. The change of every polynomial term between these points is $o\left(R_m^d\right)$, and the two remainders are also $o\left(R_m^d\right)$. Subtraction would give $\Delta E_d^+(R_m) = o\left(R_m^d\right)$, contradicting (6.8). □

Thus the obstruction is exhibited for the cube. One orthotope is sufficient to rule out a deterministic expansion universal over all side vectors; the argument does not assert the same jump multiplicities for every irrational orthotope. This is an obstruction to the aggregate expansion itself, not a gap in a discrepancy estimate. Edge memory must be exposed by exact cancellation.

### *6.4. Boolean mixed-boundary tomography*

For each coordinate, impose one boundary condition on its pair of opposite faces. Thus $\eta \in \{D, N\}^d$ denotes a separable boundary pattern, not a Zaremba condition varying within a face. Put

$$I_D = \mathbb{Z}_{>0}, \qquad I_N = \mathbb{Z}_{\geq 0}, \qquad I_\eta = \prod_{i=1}^{d} I_{\eta_i},$$

and

$$\mathcal{S}_\eta(F) = \sum_{n \in I_\eta} F\,(n).$$

For $J \subseteq [d]$, let $\eta(E)_i = N$ for $i \in E$ and $D$ otherwise, and set

$$\mathcal{T}_J(F) = \sum_{E \subseteq J} (-1)^{|J|-|E|}\, \mathcal{S}_{\eta(E)}(F). \qquad (6.9)$$

**Theorem 6.5 (exact coordinate-stratum isolation).** For every $J \subseteq [d]$, the exact isolated-stratum identity is

$$\mathcal{T}_J(F) = \sum_{\substack{n_j = 0\ (j \in J) \\ n_i \geq 1\ (i \notin J)}} F\,(n). \qquad (6.10)$$

**Proof.** In one coordinate,

$$\sum_{n_i \geq 0} f\,(n_i) - \sum_{n_i \geq 1} f\,(n_i) = f(0).$$

The differences in distinct coordinates commute. Applying them successively forces precisely the coordinates in $J$ to vanish, and expanding their product gives (6.9). □

Identity (6.10) is exact at every finite cutoff; no asymptotic remainder is introduced by the isolation step.

Figure 1. Boolean boundary-condition tomography in dimension three
The separable Dirichlet-Neumann boundary cube
DNN
NNN
DND
NND
DDN
NDN
DDD
NDD
$\Delta_2 = N_2 - D_2$
$\Delta_3 = N_3 - D_3$
$\Delta_1 = N_1 - D_1$
Each D/N choice applies to a pair of opposite coordinate faces.
Exact Boolean isolation
$\mathcal{T}_{\{i\}}(F) = \sum_{n_i=0,\ n_k\ge 1\ (k\ne i)} F(n)$
$\mathcal{T}_{\{i,j\}}(F) = \sum_{n_i=n_j=0,\ n_k\ge 1\ (k\notin\{i,j\})} F(n)$
$\mathcal{T}_{\{1,2,3\}}(F) = F(0,0,0)$
$|J| = 1$: face stratum
$|J| = 2$: edge stratum
$|J| = 3$: corner stratum

*Figure 1. The eight separable Dirichlet-Neumann patterns in dimension three. The difference $N_i - D_i$ selects $n_i = 0$, and double differences isolate codimension-two strata. The boxes use the exact transforms $\mathcal{T}_J$ of Theorem 6.5; the corner transform equals $F(0,0,0)$, which vanishes for the positive-energy weights used here. The figure is illustrative; Theorem 6.5 is the proof.*

The boundary-condition cube and the coordinate-stratum Boolean lattice are therefore the same incidence algebra. A first difference isolates a face, a double difference an edge, and a $k$-fold difference a codimension-$k$ stratum.

### 6.5. Genuine raw edge coefficients and reconstruction

Apply Theorem 6.5 to

$$F_{R,s,\Phi}(n) = c_s M_{n,a}^s \Phi\big(q(n;a)\big)\mathbf{1}_{\{0<M_{n,a}\le R\}}.$$

For $J \subsetneq [d]$, put

$$m = d - |J|, \qquad V_{J^c}(a) = \prod_{i\notin J} a_i,$$

and let $\iota_J$ insert zero coordinates on $J$.

**Theorem 6.6 (recursive stratum-resolved raw measure).** If $m \ge 2$, then

$$\begin{aligned}\mathcal{T}_J\big(F_{R,s,\Phi}\big) = \frac{c_s V_{J^c}(a)|S^{m-1}|}{2^m(m+s)} R^{m+s}\int_{\Delta_{m-1}} \Phi\left(\iota_J x\right) d\mu_m(x) \\ -\frac{c_s V_{J^c}(a)|S^{m-2}|}{2^m(m+s-1)} R^{m+s-1}\sum_{i\notin J}\frac{1}{a_i}\int_{\Delta_{m-2}} \Phi\left(\iota_{J\cup\{i\}}x\right) d\mu_{m-1}(x) \\ +O_{a,\Phi}\big(R^{m+s-2+2/(m+1)}\big).\end{aligned} \tag{6.11}$$

For $m = 1$, the same leading term holds with remainder $O_{a,\Phi}(R^s)$.

**Proof.** Theorem 6.5 reduces the transform to the strictly positive orthant in the $m$ surviving coordinates. Apply Lemma 2.4 in precisely those coordinates. Its full-lattice term gives the first line of (6.11). Its codimension-one terms have negative sign and give the second line, with lattice density $V_{J^c}(a)/a_i$ and angular law $\mu_{m-1}$ on the inserted face. Lemma 6.1 bounds the full-stratum remainder by the final line of (6.11). The face remainders and every stratum of codimension at least two have smaller order and are absorbed there. Direct summation treats $m = 1$. □

For $s = 2, J = \{j, k\}$, and $d \geq 4$, (6.11) becomes

$$\begin{aligned} \mathcal{T}_{\{j,k\}}\left(F_{R,2,\Phi}\right) = {} & \frac{\pi^2 V_{\{j,k\}^c}(a)\left|S^{d-3}\right|}{2^{d-2}d} R^d \int_{\Delta_{d-3}} \Phi\left(\iota_{\{j,k\}}x\right) d\mu_{d-2}(x) \\ & - \frac{\pi^2 V_{\{j,k\}^c}(a)\left|S^{d-4}\right|}{2^{d-2}(d-1)} R^{d-1} \sum_{i\notin\{j,k\}} \frac{1}{a_i} \int_{\Delta_{d-4}} \Phi\left(\iota_{\{i,j,k\}}x\right) d\mu_{d-3}(x) \\ & + O_{a,\Phi}\left(R^{d-2+2/(d-1)}\right). \end{aligned} \tag{6.12}$$

For $d = 3$, the surviving dimension is $m = 1$, so the first line of (6.12) holds with remainder $O_{a,\Phi}(R^2)$ and there is no nonzero lower stratum. Thus the $R^d$ term is a genuine raw edge coefficient after exact cancellation, although Theorem 6.4 forbids it as a third coefficient of the uncancelled population. For $d \geq 4$, the additional $R^{d-1}$ term is precisely the face coefficient of the isolated edge stratum. The same Boolean geometry therefore recurs inside every isolated coordinate stratum.

The scalar leading coefficients are

$$b_{jk} = V_{\{j,k\}^c}(a) = \frac{V(a)}{a_j a_k}.$$

Their collection reconstructs the orthotope:

$$V(a) = \left(\prod_{j<k} b_{jk}\right)^{2/((d-1)(d-2))}, \qquad d \geq 3, \tag{6.13}$$

and, for distinct $i, j, k$,

$$a_i = \left(\frac{V(a) b_{jk}}{b_{ij} b_{ik}}\right)^{1/2}. \tag{6.14}$$

Equations (6.13) and (6.14) follow by substituting $b_{jk} = V/(a_j a_k)$. This reconstruction uses the mixed-boundary data entering the Boolean transforms, namely the four boundary patterns associated with each pair $\{j, k\}$. The raw mixed-boundary family includes the pure Dirichlet spectrum and therefore includes the input from which the heat invariants of Proposition 4.3 are extracted. The finite coefficient systems are nevertheless different: Proposition 4.3 uses $d$ coefficients of one complete Dirichlet heat trace, whereas (6.13)–(6.14) use one isolated asymptotic coefficient for each coordinate pair across mixed boundary patterns. Each system reconstructs the same side vector and hence determines the other through that reconstruction; their distinction is the extraction mechanism and the geometric stratum made explicit. The present result shows that the obstruction in Theorem 6.4 is caused by aggregation rather than loss of geometry, and that the coordinate-stratum hierarchy remains resolvable after exact boundary-condition inversion.

### *6.6. Transform stability and the exact boundary*

For either $\eta = D$ or $\eta = N$, define the first Riesz and heat transforms of the raw tensor by

$$\mathsf{R}_a^{\eta}(\Lambda) = \int_0^{\Lambda} \mathsf{E}_a^{\eta}(u)\,du, \qquad \mathsf{H}_a^{\eta}(t) = t\int_0^{\infty} e^{-tu}\,\mathsf{E}_a^{\eta}(u)\,du.$$

Termwise integration of (6.6) gives

$$\mathsf{R}_a^{D/N}(\Lambda) = \frac{C_d}{d/2+2} V(a)\Lambda^{d/2+2} I$$
$$\mp \frac{D_d}{d+3}\mathsf{F}(K_a)\Lambda^{(d+3)/2} + O_a\big(\Lambda^{d/2+1+1/(d+1)}\big).$$

Here $D_d/(d+3) = (D_d/2)/\big((d+3)/2\big)$. Laplace integration gives, as $t \downarrow 0$,

$$\mathsf{H}_a^{D/N}(t) = C_d\Gamma\left(\frac{d}{2}+2\right)V(a)t^{-d/2-1} I$$
$$\mp \frac{D_d}{2}\Gamma\left(\frac{d+3}{2}\right)\mathsf{F}(K_a)t^{-(d+1)/2} + O_a\big(t^{-d/2-1/(d+1)}\big).$$

These formulas preserve the same face-memory tensor. They are consequences of the raw theorem, not substitutes for it; general scalar Riesz-mean asymptotics provide the surrounding context [8].

The absolute-summability extension in Lemma 2.4 is used here without truncation. For

$$F_{t,s,\Phi}(n) = c_s M_{n,a}^s \Phi\big(q(n;a)\big) e^{-\pi^2 t M_{n,a}^2}, \qquad t > 0,$$

the full lattice series is absolutely convergent. Lemma 2.4 decomposes its Dirichlet sum into all coordinate strata, while Theorem 6.5 isolates any one of them by mixed boundary differences. Thus the raw, Riesz, and heat populations are three transforms of the same Boolean stratum calculus.

The section resolves universality and memory at consecutive levels. The bulk gives the universal simplex law. The singular face part of the bulk–face correction, and its tensor moment, retain the complete side geometry. Edge-sized jumps prevent an uncancelled third raw coefficient, but exact mixed-boundary Möbius inversion isolates every coordinate stratum and produces a genuine edge coefficient. Within the isolated stratum, the next lower-dimensional coefficient and a remainder below that scale are again explicit.

This cancellation mechanism is specific to separable orthotopes. Section 8 replaces the coordinate Boolean lattice by the incidence lattice of a great-circle arrangement. There the relevant aggregation is not a boundary-condition transform, but the signed superposition of local spherical pushforwards.

# 7. Identifiability beyond orthotopes

Section 6 recovers an orthotope from a tensorial coefficient of its directional spectral projector. We now consider a different inverse problem: recovery from direction-labelled width data. The distinction is essential. The results below do not identify a general convex body from an unlabelled width distribution, and they do not assert that the required directional function is encoded by an ordinary unlabelled spectrum.

### 7.1. The exact information boundary

For a convex body $K \subset \mathbb{R}^d$, let $h_K$ be its support function and let

$$D_K = \frac{1}{2}(K - K)$$

be its central symmetral. Additivity of support functions gives

$$w_K(u) = h_K(u) + h_K(-u) = 2h_{D_K}(u).$$

Consequently, the complete direction-labelled width function determines $D_K$, but cannot distinguish bodies with the same difference body. If $K$ is centrally symmetric, it determines $K$ up to translation. Passing from the function $u \mapsto w_K(u)$ to the law of $w_K(U)$ forgets the assignment of values to directions and introduces a further obstruction.

That second loss occurs even among smooth, strictly convex, centrally symmetric plane bodies. For $0 < \varepsilon < 1/15$, consider the support functions

$$h_1(\theta) = 1 + \varepsilon\cos 2\theta, \qquad h_2(\theta) = 1 + \varepsilon\cos 4\theta.$$

Since

$$h_1 + {h_1}'' = 1 - 3\varepsilon\cos 2\theta > 0, \qquad h_2 + {h_2}'' = 1 - 15\varepsilon\cos 4\theta > 0,$$

the curves

$$x_k(\theta) = h_k(\theta)n(\theta) + {h_k}'(\theta)t(\theta), \qquad n = (\cos\theta, \sin\theta), \quad t = (-\sin\theta, \cos\theta),$$

satisfy ${x_k}' = (h_k + {h_k}'')t$. They therefore turn strictly monotonically and bound smooth strictly convex bodies having $h_k$ as support functions. The identity $h_k(\theta + \pi) = h_k(\theta)$ gives central symmetry. If $\Theta$ is uniform on $[0,2\pi)$, then $2\Theta$ and $4\Theta$, reduced modulo $2\pi$, are both uniform, and hence

$$\mathcal{L}\big(2h_1(\Theta)\big) = \mathcal{L}\big(2h_2(\Theta)\big).$$

The two bodies are not congruent. A translation changes a support function only by a first spherical harmonic, while a planar rotation only changes the phase of a fixed Fourier order; neither operation can transform a nonzero second harmonic into a fourth harmonic. Thus unlabelled width laws are not injective even on that restricted class. This agrees with the broader distinction between labelled and distributional projection data in geometric tomography; see [9, 16, 21].

The remainder of the section gives positive theorems of two kinds. First, the complete width law is injective on certain finite-dimensional classes. Second, the complete direction-labelled function admits an exact canonical decomposition on the full linear–quadratic class.

### 7.2. The canonical linear–quadratic class

Let $g_1, \dots, g_N \in \mathbb{R}^d \setminus \mathbf{0}$, and let $A$ be symmetric and positive semidefinite. Define

$$\mathcal{W}_{G,A}(u) = \sum_{j=1}^{N} |g_j \cdot u| + 2\sqrt{u^T A u}, \qquad u \in S^{d-1}. \tag{7.1}$$

If

$$Z_G = \sum_{j=1}^{N} [-g_j/2, g_j/2], \qquad E_A = A^{1/2} B_2^d,$$

then $\mathcal{W}_{G,A}$ is the width function of $Z_G + E_A$. Indeed,

$$h_{Z_G}(u) = \frac{1}{2} \sum_j |g_j \cdot u|, \qquad h_{E_A}(u) = \sqrt{u^T A u},$$

and support functions add under Minkowski addition.

Parallel generators are not individually identifiable, since their absolute linear forms add. We therefore combine every parallel class into a single generator and call the resulting family *reduced*. There is one further, intrinsic ambiguity. If $A = \alpha v \otimes v$ has rank one, then

$$2\sqrt{u^T A u} = |2\sqrt{\alpha}\, v \cdot u|,$$

so the quadratic term is itself a linear generator. We absorb such a term into the reduced generator family. A representation is called *canonical* if the generator family is reduced and either $A = 0$ or $\operatorname{rank} A \geq 2$. These reductions make the representation in (7.1) intrinsic.

For $A = 0$, uniqueness of the reduced generator measure is also an instance of the classical injectivity of the spherical cosine transform on even measures, equivalently uniqueness of the generating measure of a zonoid [21]. The proof below is retained because it separates the atomic linear part constructively from the quadratic residual by ridge locations and gradient jumps. That constructive separation, and its finite moment form in Theorem 7.3, are the mechanisms used later.

The same zonoid interpretation extends to the full quadratic term. Let $\sigma$ be normalized spherical measure and put

$$c_d = \left( \int_{S^{d-1}} |e_1 \cdot v| \, d\sigma(v) \right)^{-1}.$$

Rotational invariance gives, for every $x \in \mathbb{R}^d$,

$$|x| = c_d \int_{S^{d-1}} |x \cdot v| \, d\sigma(v).$$

Applying this identity to $x = A^{1/2} u$ shows that $E_A$ is a zonoid:

$$h_{E_A}(u) = c_d \int_{S^{d-1}} |u \cdot A^{1/2} v| \, d\sigma(v) = \int_{S^{d-1}} |u \cdot w| \, d\nu_A(w), \tag{7.2}$$

where $\nu_A$ is the even weighted pushforward defined by

$$\int f(w)\, d\nu_A(w) = c_d \int_{\{A^{1/2}v \neq 0\}} \left|A^{1/2}v\right| f\left(\frac{A^{1/2}v}{\left|A^{1/2}v\right|}\right) d\sigma(v).$$

If $\mathrm{rank}A \geq 2$, the inverse image of a fixed direction under $v \mapsto A^{1/2}v/\left|A^{1/2}v\right|$ has spherical measure zero; hence $\nu_A$ is atomless. If $\mathrm{rank}A = 1$, the image consists of one antipodal pair and $\nu_A$ is atomic, exactly as the segment identity above predicts. The finite zonotopal part has an atomic generating measure, while the canonical quadratic part has an atomless one. By uniqueness of the even generating measure under the cosine transform [21], its Lebesgue decomposition therefore separates these two parts uniquely. Thus canonicality is precisely the condition under which the classical zonoid decomposition and the constructive ridge decomposition below agree.

**Theorem 7.1 (canonical linear–quadratic identifiability).** The complete direction-labelled function $\mathcal{W}_{G,A}$ determines every canonical pair $(G, A)$ uniquely, up to permutation and sign changes of the generators. Equivalently, it determines the centred body $Z_G + E_A$; without a prescribed centre it determines the body up to translation.

**Proof.** For each generator put

$$H_j = g_j^{\perp} \cap S^{d-1}.$$

The absolute-linear part is smooth away from $\bigcup_j H_j$. The quadratic part is smooth away from

$$\Sigma_A = S^{d-1} \cap \ker A.$$

If $A \neq 0$, the codimension of $\Sigma_A$ in the sphere is $\mathrm{rank}A$. Canonicality makes this codimension at least two. It follows that the codimension-one components of the nondifferentiability set of $\mathcal{W}_{G,A}$ are precisely the distinct great spheres $H_j$. Indeed, $H_j \subseteq \Sigma_A$ would imply $g_j^{\perp} \subseteq \ker A$, hence $\mathrm{rank}A \leq 1$, contrary to canonicality.

Fix a point $u \in H_j$ outside every other $H_k$ and outside $\Sigma_A$. Such points form a dense open subset of $H_j$. Let $\nabla_S^+$ and $\nabla_S^-$ denote the one-sided spherical gradients from the two sides of $H_j$. Every term except $\left|g_j \cdot u\right|$ has the same limiting gradient from both sides. Since $g_j \cdot u = 0$, the tangent projection of $g_j$ at $u$ is $g_j$ itself, and therefore

$$\nabla_S^+ \mathcal{W}_{G,A}(u) - \nabla_S^- \mathcal{W}_{G,A}(u) = 2g_j$$

after orienting the two sides. Reversing that orientation changes only the sign. Thus the ridge location recovers the line $\mathbb{R}g_j$, and the gradient jump recovers $\left|g_j\right|$. The entire reduced generator multiset $\{\pm g_1, \dots, \pm g_N\}$ is consequently determined.

Subtract its recovered contribution from (7.1). The residual is

$$r_A(u) = 2\sqrt{u^T A u}.$$

It determines the quadratic form on all of $\mathbb{R}^d$ by positive homogeneity:

$$Q_A(x) = \begin{cases} \dfrac{|x|^2}{4} r_A\left(\dfrac{x}{|x|}\right)^2, & x \neq 0, \\ 0, & x = 0. \end{cases}$$

Finally,

$$x^T Ay = \frac{1}{4}\big(Q_A(x+y) - Q_A(x-y)\big)$$

recovers every matrix entry of $A$. Hence both parts of the canonical decomposition are unique. □

The rank condition in Theorem 7.1 is sharp as a statement about the decomposition: rank-one quadratic terms and single linear generators are the same functions. It is not a technical exclusion.

**Corollary 7.2.** Theorem 7.1 includes the following direction-labelled classes in one recovery theorem:

1. reduced zonotopes, by taking $A = 0$;
2. ellipsoids, by taking $N = 0$;
3. right circular cylinders in $\mathbb{R}^3$, by taking $g_1 = He$ and $A = (D^2/4)(I - e \otimes e)$;
4. arbitrary Minkowski sums of a reduced zonotope and an ellipsoid whose quadratic matrix has rank at least two.

For orthotopes, the generator directions recovered here are the same principal directions in which the face-memory tensor in (6.6) is diagonal. The two theorems use different data and neither is substituted for the other: Section 6 recovers side geometry from a spectral-projector coefficient, while Theorem 7.1 recovers a canonical mixed body from its direction-labelled width function.

### *7.3. Finite direction-sensitive identifiability*

Theorem 7.1 uses the complete direction-labelled function. Under a fixed bound on generator complexity, a finite set of distributional ridge moments already suffices. These moments retain directional information and must not be confused with scalar moments of the random variable $\mathcal{W}_{G,A}(U)$.

Let $g_j = \ell_j v_j$, where $\ell_j > 0$, $|v_j| = 1$, and put $P_j = v_j \otimes v_j$. Denote by $\mathcal{S}_{\mathcal{W}}$ the hypersurface-singular part of $(\Delta_S + d - 1)\mathcal{W}_{G,A}$.

**Theorem 7.3 (finite ridge-moment identifiability).** Suppose that a canonical pair has at most $N_0$ reduced generators. The tensors

$$\mathsf{H}_k(\mathcal{W}) = \int_{S^{d-1}} u^{\otimes 2k}\, d\mathcal{S}_{\mathcal{W}}(u), \qquad 0 \le k \le 2N_0 - 1, \qquad (7.3)$$

together with one residual quadratic tensor moment, determine the complete canonical pair up to permutation and sign changes of its generators.

**Proof.** Integration by parts over the two hemispheres separated by $v_j^{\perp}$ gives

$$(\Delta_S + d - 1)|g_j \cdot u| = 2\ell_j\, \mathcal{H}^{d-2} \restriction \left(v_j^{\perp} \cap S^{d-1}\right)$$

in its hypersurface-singular part. The quadratic term has no codimension-one singular measure: its singular set is empty when $A$ is positive definite and, by the canonical rank condition, has codimension at least two otherwise. Hence

$$\mathcal{S}_{\mathcal{W}} = 2\sum_j \ell_j \, \mathcal{H}^{d-2} \restriction \left(v_j^{\perp} \cap S^{d-1}\right). \qquad (7.4)$$

Equation (7.4) separates the codimension-one singular measure from the quadratic residual.

Rotational symmetry on $v_j^{\perp}$ gives, for every $x \in \mathbb{R}^d$,

$$\int_{v_j^{\perp} \cap S^{d-1}} (u \cdot x)^{2k} \, d\mathcal{H}^{d-2}(u) = c_{d,k}(x^T\left(I - P_j\right)x)^k.$$

Only the nonvanishing of the dimensional constant $c_{d,k}$ is used.

Expansion of the right-hand side is triangular in $\left(x^T P_j x\right)^0, \dots, \left(x^T P_j x\right)^k$, with nonzero top coefficient. Therefore the tensors in (7.3) determine

$$\mathsf{J}_k = \sum_j \ell_j \, P_j^{\otimes k}, \qquad 0 \le k \le 2N_0 - 1. \qquad (7.5)$$

We prove finite determinacy directly. Suppose two atomic measures, each with at most $N_0$ atoms, have the moments (7.5). Their difference has at most $2N_0$ support points $x_1, \dots, x_s$. For fixed $i$, choose affine functions $L_{ij}$ satisfying

$$L_{ij}(x_i) = 1, \qquad L_{ij}\left(x_j\right) = 0,$$

and put

$$p_i = \prod_{j \neq i} L_{ij}.$$

Its degree is at most $2N_0 - 1$. Equality of the truncated moments makes its integral against the difference measure vanish, but this integral is exactly the weight at $x_i$. Every weight therefore vanishes. Applied to (7.5), this recovers every projector $P_j$ and every length $\ell_j$, hence every generator up to sign and permutation.

Subtract the recovered linear part and write

$$r(u) = 2\sqrt{u^T A u}, \qquad \mathsf{Q} = \int_{S^{d-1}} r\,(u)^2 u \otimes u \, d\sigma(u).$$

Here $\sigma$ is normalized spherical measure.

The fourth spherical moment identity gives

$$\mathsf{Q} = \frac{4}{d(d+2)}\left((\mathrm{tr} A)I + 2A\right).$$

Taking the trace and solving for $A$ yields

$$A = \frac{d}{8}\left((d+2)\mathsf{Q} - (\mathrm{tr}\mathsf{Q})I\right). \qquad (7.6)$$

Formula (7.6) completes the recovery of the quadratic part. Thus the finite data determine the complete canonical pair. □

The order in Theorem 7.3 depends on the generator bound $N_0$. The theorem does not assert generator-independent scalar-moment closure, injectivity of the unlabelled scalar width law on the full mixed class, or reconstruction from an ordinary Laplacian spectrum.

### *7.4. Law-level recovery for cylinders and ellipsoids*

The obstruction of Subsection 7.1 concerns unrestricted classes. On the following finite-dimensional families, even the unlabelled width law is injective.

**Proposition 7.4 (right circular cylinders).** Let $C_{H,D} \subset \mathbb{R}^3$ be a centred right circular cylinder of height $H > 0$, base diameter $D > 0$, and axis $e$. Put

$$L = \sqrt{H^2 + D^2}, \qquad m = \min(H, D), \qquad M = \max(H, D).$$

Its random width has the representation

$$W_{H,D} = HT + D\sqrt{1 - T^2}, \qquad T \sim \text{Unif}[0,1].$$

Writing $q(w) = \sqrt{L^2 - w^2}$, its density vanishes outside $[m, L]$. If $H < D$, then

$$f_{H,D}(w) = \begin{cases} \dfrac{Dw/q(w) - H}{L^2}, & H < w < D, \\ \dfrac{2Dw}{L^2 q(w)}, & D < w < L, \end{cases}$$

whereas, if $H > D$,

$$f_{H,D}(w) = \begin{cases} \dfrac{H + Dw/q(w)}{L^2}, & D < w < H, \\ \dfrac{2Dw}{L^2 q(w)}, & H < w < L. \end{cases}$$

When $H = D$, only the second branch occurs on $(H, L)$. The law uniquely determines the ordered pair $(H, D)$.

**Proof.** The cylinder is the Minkowski sum of an axial segment of length $H$ and a radius-$D/2$ disk orthogonal to $e$. Its width is therefore

$$H|U \cdot e| + D\sqrt{1 - (U \cdot e)^2}.$$

For uniform $U \in S^2$, the variable $T = |U \cdot e|$ is uniform on $[0,1]$. The function $g(t) = Ht + D\sqrt{1 - t^2}$ has endpoint values $D, H$, a unique maximum $L$ at $t = H/L$, and inverse branches

$$t_\pm(w) = \frac{Hw \pm Dq(w)}{L^2}, \qquad t_\pm{}'(w) = \frac{H \mp Dw/q(w)}{L^2}.$$

On $(m, M)$, exactly one branch satisfies the unsquared defining equation; on $(M, L)$, both do. Summing their absolute derivatives proves the density formulae.

The support supplies $m$ and $L$, hence $M = \sqrt{L^2 - m^2}$. If $H < D$, the density tends to zero at the lower endpoint. If $H > D$, it tends to $1/M$. The equality case is characterized by $L^2 = 2m^2$. Thus the law distinguishes which of $m, M$ is the height and which is the diameter. □

At $w = L$, the density has the fold

$$f_{H,D}(w) \sim \frac{\sqrt{2}D}{L^{3/2}}(L - w)^{-1/2}.$$

The internal transition is a positive step when $H < D$ and a continuous corner when $H > D$. The singularity type is therefore part of the inversion mechanism, not an incidental feature of the formula.

**Proposition 7.5 (ellipsoid moment reconstruction).** Let $\mathcal{E}_a \subset \mathbb{R}^d$ be the centred ellipsoid with positive semiaxes $a_1, \dots, a_d$. Its first $d$ even width moments determine the unordered semiaxis vector.

**Proof.** The width is

$$W_a(U) = 2\left(\sum_{i=1}^{d} a_i^2 U_i^2\right)^{1/2}.$$

Put $M_m = \mathbb{E}W_a^{2m}$, $M_0 = 1$, and

$$c_m = \frac{(d/2)_m}{4^m m!} M_m.$$

By Proposition 2.1 and the multinomial theorem,

$$\sum_{m\geq 0} c_m t^m = \prod_{i=1}^{d} \left(1 - a_i^2 t\right)^{-1/2}. \qquad (7.7)$$

Taking the formal logarithm of (7.7) gives

$$\log\left(\sum_{m\geq 0} c_m t^m\right) = \frac{1}{2}\sum_{k\geq 1} \frac{t^k}{k} \sum_{i=1}^{d} a_i^{2k}.$$

The first $d$ even moments therefore determine the first $d$ power sums of $a_1^2, \dots, a_d^2$. Newton's identities recover their elementary symmetric functions and hence the monic polynomial having those numbers as roots. Positivity then recovers the semiaxes. □

For a three-dimensional spheroid with semiaxes $(R, R, c)$, put $T = |U_3|$. Then $T$ is uniform on $[0,1]$ and

$$W = 2\sqrt{R^2 + (c^2 - R^2)T^2}.$$

Solving this monotone relation for $T$ and differentiating gives the following explicit laws. If $c > R$, the density on $(2R, 2c)$ is

$$f(w) = \frac{w}{2\sqrt{c^2 - R^2}\sqrt{w^2 - 4R^2}},$$

whereas, if $c < R$, the density on $(2c, 2R)$ is

$$f(w) = \frac{w}{2\sqrt{R^2 - c^2}\sqrt{4R^2 - w^2}}.$$

The spherical case is the atom at $2R$. The support endpoints and the endpoint carrying the inverse-square-root fold recover $R$ and $c$.

### *7.5. A common moment operator*

The linear and quadratic members also share one exact radialization formula. Let $G_0 \sim N(0, I_d)$, write $G_0 = RU$, and let $p \geq 0$ be an integer. Homogeneity and independence of $R$ and $U$ give

$$\mathbb{E}\,\mathcal{W}_{G,A}(U)^p = \frac{\Gamma(d/2)}{2^{p/2}\Gamma\big((d+p)/2\big)}\mathbb{E}\left(\sum_j |g_j \cdot G_0| + 2\sqrt{G_0^T A G_0}\right)^p. \qquad (7.8)$$

Indeed, the expression inside the last expectation is $R\mathcal{W}_{G,A}(U)$, and

$$\mathbb{E}R^p = 2^{p/2}\frac{\Gamma\big((d+p)/2\big)}{\Gamma(d/2)}.$$

Equation (7.8) extends the Gaussian radialization of Section 3 to the complete linear–quadratic class. It does not factor the correlated Gaussian terms and is not used to infer identifiability. The latter comes from the ridge-jump and polarization argument of Theorem 7.1. Together, Theorems 7.1 and 7.3 and Propositions 7.4–7.5 locate the precise positive boundary: canonical mixed bodies are identifiable from direction-labelled widths, while certain distinguished families are already identifiable from their unlabelled laws; neither statement extends to unrestricted convex bodies.

# 8. Stratified zonotopal width laws

Section 6 recovered an orthotope from coordinate-stratum coefficients of its spectral projector, while Section 7 recovered a canonical linear–quadratic body from its complete direction-labelled width function. We now pass from the labelled function to its scalar pushforward law and determine exactly which features of a reduced zonotopal arrangement remain visible after directions have been forgotten.

For orthotopes, the strata form a coordinate Boolean lattice and can be isolated by mixed Dirichlet-Neumann differences. For a general zonotope, the coordinate lattice is replaced by the incidence lattice of a great-circle arrangement. There is no analogue here of the boundary-condition inversion of Section 6. The common principle is that the geometric strata must be resolved before their contributions are aggregated. Once aggregation has occurred, distinct local contributions at the same width value may reinforce or cancel.

The support-function representation of zonotopes is classical [21], as are co-area [6], hyperplane-arrangement language [20], and the general stratified-critical-point viewpoint [11]. The theorem proved below is the explicit synthesis required here: a global probability density, complete incidence tests, one common local normal form, normalized step, fold, endpoint-fold and corner coefficients, resonance rules, signed superposition, and an exact cancellation counterexample.

### *8.1. Reduced generators and global co-area*

Let

$$Z = \sum_{k=1}^{N}[-g_k/2, g_k/2] \subset \mathbb{R}^3$$

be full-dimensional. Parallel generators are first oriented consistently and merged. The resulting family is pairwise nonparallel and is called reduced. For $U$ uniform on $S^2$, put

$$W(U) = w_Z(U) = \sum_{k=1}^{N} |g_k \cdot U|.$$

Let

$$\mathcal{A} = \{g_k^\perp \cap S^2 : 1 \le k \le N\}$$

be the great-circle arrangement. An open two-cell $C$ has a sign vector $\varepsilon(C) = (\varepsilon_1, \dots, \varepsilon_N)$, and we write

$$G_C = \sum_{k=1}^{N} \varepsilon_k \, g_k.$$

Then $W(u) = G_C \cdot u$ on $C$. Since $W > 0$ on $S^2$, no nonempty two-cell has $G_C = 0$.

For $0 < \rho < |G_C|$, the level set $G_C \cdot u = \rho$ is a Euclidean circle of radius

$$r_C(\rho) = \frac{\sqrt{|G_C|^2 - \rho^2}}{|G_C|}.$$

The spherical gradient of $G_C \cdot u$ has length $\sqrt{|G_C|^2 - \rho^2}$ on that circle. Arc length is $r_C(\rho)\, d\theta$; hence the co-area quotient is exactly $d\theta / |G_C|$. If $\Psi_C(\rho)$ is the angular measure of the part of the circle in $C$, summation over the cells gives

$$f(\rho) = \frac{1}{4\pi} \sum_C \frac{\Psi_C(\rho)}{|G_C|}. \qquad (8.1)$$

Cell boundaries have spherical area zero, so this also proves absolute continuity. Formula (8.1) is derived from co-area and is not a regularity assumption.

Integrating the cellwise co-area formula recovers the spherical probability measure exactly:

$$\int_{\mathbb{R}} f(\rho)\, d\rho = \frac{1}{4\pi} \sum_C \operatorname{area}(C) = \frac{1}{4\pi} \operatorname{area}(S^2) = 1.$$

For an orthogonal box, the eight cell gradients have the common norm $d = (a_1^2 + a_2^2 + a_3^2)^{1/2}$, and their angular contributions agree by reflection. Hence (8.1) reduces to

$$f(\rho) = \frac{2}{\pi d} \Psi(\rho),$$

which is the global co-area density of the companion paper. Thus (8.1) is a strict arrangement extension of the earlier box formula.

Figure 2. A reduced great-circle arrangement and its incidence strata

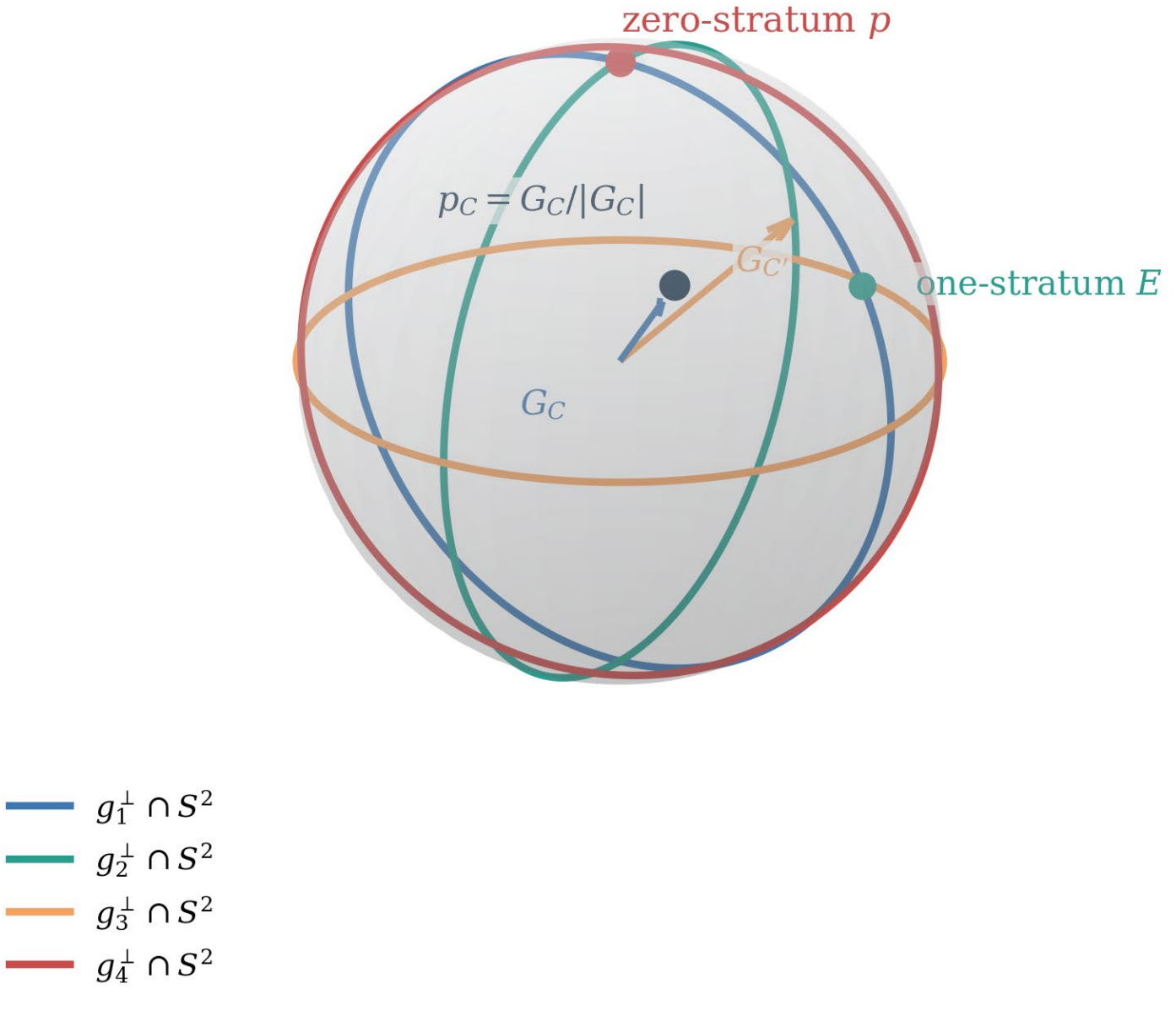

On each sign cell $C$, the zonotopal width is the spherical linear form $W_Z(u) = G_C \cdot u$.

*Figure 2. A reduced great-circle arrangement, its sign cells, one-strata, zero-strata, and adjacent cell gradients. The figure illustrates the objects in the global co-area formula; incidence is decided algebraically below.*

### 8.2. Complete incidence on every stratum

Fix a two-cell $C$. The positive stationary point of its linear restriction to the sphere is

$$p_C = \frac{G_C}{|G_C|}, \qquad c_C = |G_C|.$$

It belongs to $C$, respectively to $\overline{C}$, precisely when

$$\varepsilon_k\, g_k \cdot G_C > 0, \qquad \text{respectively} \qquad \varepsilon_k\, g_k \cdot G_C \geq 0 \quad (1 \leq k \leq N). \tag{8.2}$$

Now let $E$ be an open one-stratum contained in $g_j^\perp \cap S^2$. The signs $\varepsilon_k$ are fixed on $E$ for $k \neq j$. Set

$$H_E = \sum_{k \neq j} \varepsilon_k\, g_k, \qquad \Pi_j = I - \frac{g_j g_j^T}{|g_j|^2}.$$

On $E$, $W(u) = H_E \cdot u = (\Pi_j H_E) \cdot u$. The vector $\Pi_j H_E$ cannot vanish: otherwise $W$ would vanish on $E$, contrary to the full dimensionality of $Z$. The only stationary point compatible with $W > 0$ is therefore

$$p_E = \frac{\Pi_j H_E}{|\Pi_j H_E|}, \qquad c_E = |\Pi_j H_E|.$$

It lies in $E$, respectively $\overline{E}$, exactly when

$$\varepsilon_k\, g_k \cdot p_E > 0, \qquad \text{respectively} \qquad \varepsilon_k\, g_k \cdot p_E \geq 0 \quad (k \neq j). \qquad (8.3)$$

Finally, a zero-stratum $p$ is the common intersection of the great circles whose active index set is

$$I(p) = \{k : g_k \cdot p = 0\}. \qquad (8.4)$$

Its attained value is $c_p = W(p)$. Equations (8.2), (8.3), and (8.4) are the complete incidence test. Equal numerical values are combined only after these tests have been applied.

The same incidence calculus gives the support without any genericity assumption. On an open two-cell, the only positive stationary point of the spherical linear form is $p_C$, and it is a maximum. On an open one-stratum, the only positive stationary point is $p_E$, again a maximum of the restriction. A global minimum must therefore occur at a zero-stratum. For the upper endpoint, write

$$W(u) = \max_{\varepsilon \in \{-1,1\}^N} \left( \sum_j \varepsilon_j\, g_j \right) \cdot u.$$

Maximizing first in $u$ gives the largest norm of a signed sum. If a maximizing sign vector were not realized at the direction of its signed sum, the realized sign vector there would give a strictly larger norm, a contradiction. Hence a maximizing sign vector satisfies the closed incidence test (8.2), and

$$\operatorname{supp}\mathcal{L}(W) = \left[ \min_{p \text{ zero-stratum}} W(p),\ \max_C |G_C| \right].$$

### *8.3. Analyticity away from the finite critical set*

For a fixed cell put $G = G_C$, $\gamma = |G|$, and $n = G/\gamma$. Choose an orthonormal basis $e_1, e_2$ of $n^\perp$. Its level circle has the parametrization

$$u(\rho, \theta) = \frac{\rho}{\gamma} n + \sqrt{1 - \frac{\rho^2}{\gamma^2}}\, (\cos\theta\, e_1 + \sin\theta\, e_2).$$

For every boundary inequality of $C$,

$$F_k(\rho, \theta) = \varepsilon_k g_k \cdot u(\rho, \theta) > 0$$

is real-analytic for $|\rho| < \gamma$. A boundary angle is consequently a real-analytic function of $\rho$ whenever $F_k = 0$ and $\partial_\theta F_k \neq 0$. Its analyticity can fail only in one of the following ways.

1. $F_k = \partial_\theta F_k = 0$. The level circle is tangent to $g_k^\perp \cap S^2$, and the contact point is the stationary point (8.3) of the corresponding one-stratum.
2. Boundary angles belonging to two or more distinct great circles collide. Their common point is a zero-stratum and $\rho = W(p)$.
3. The level circle shrinks at $\rho = \gamma$. It meets $\overline{C}$ in the limit precisely under (8.2).

On a compact $\rho$-interval avoiding those values, every boundary angle is analytic, its cyclic order is fixed, and $\Psi_C$ is a finite sum of differences of analytic endpoint angles. Thus $\Psi_C$, and hence $f$, is real-analytic there. There are finitely many cells, edges, and vertices, so the exceptional set is finite.

### *8.4. The common local normal form*

Fix $p \in S^2$ and identify $T_pS^2$ with $\mathbb{R}^2$. Gnomonic coordinates are

$$u(x) = \frac{p + x}{\sqrt{1 + |x|^2}}. \qquad (8.5)$$

Their surface element is

$$dS\big(u(x)\big) = (1 + |x|^2)^{-3/2}\, dx.$$

Indeed, in an orthonormal tangent basis the Gram determinant of the two derivatives of (8.5) is $(1 + |x|^2)^{-3}$.

If $C$ is incident with $p$, put $c = W(p) = G_C \cdot p$ and $q_C = G_C - cp \in T_pS^2$. On that cell,

$$W\big(u(x)\big) - c = \frac{c + q_C \cdot x}{\sqrt{1 + |x|^2}} - c = q_C \cdot x - \frac{c}{2}|x|^2 + O(|x|^3). \qquad (8.6)$$

If $k \in I(p)$, then $g_k \cdot u(x)$ has the sign of $g_k \cdot x$. Thus all incident cell boundaries are exactly straight rays in these coordinates, not merely tangent approximations. Equivalently, before a cell is selected,

$$W\big(u(x)\big) - c = \ell_p \cdot x + \sum_{k \in I(p)} |g_k \cdot x| - \frac{c}{2}|x|^2 + O(|x|^3), \qquad (8.7)$$

where

$$\ell_p = \sum_{k \notin I(p)} \operatorname{sgn}(p \cdot g_k)(g_k - (p \cdot g_k)p).$$

Formula (8.7) is the common local engine. If the stratum through $p$ has dimension $s$, its $2 - s$ normal variables occur linearly, while its $s$ tangential variables occur quadratically at a critical point. A value window of size $\delta$ therefore has area order

$$\delta^{\,2-s}\delta^{\,s/2} = \delta^{\,2-s/2},$$

and density order $\delta^{\,1-s/2}$. The next three sections compute the coefficients and all degeneracies, so this scaling calculation is not used as a substitute for proof.

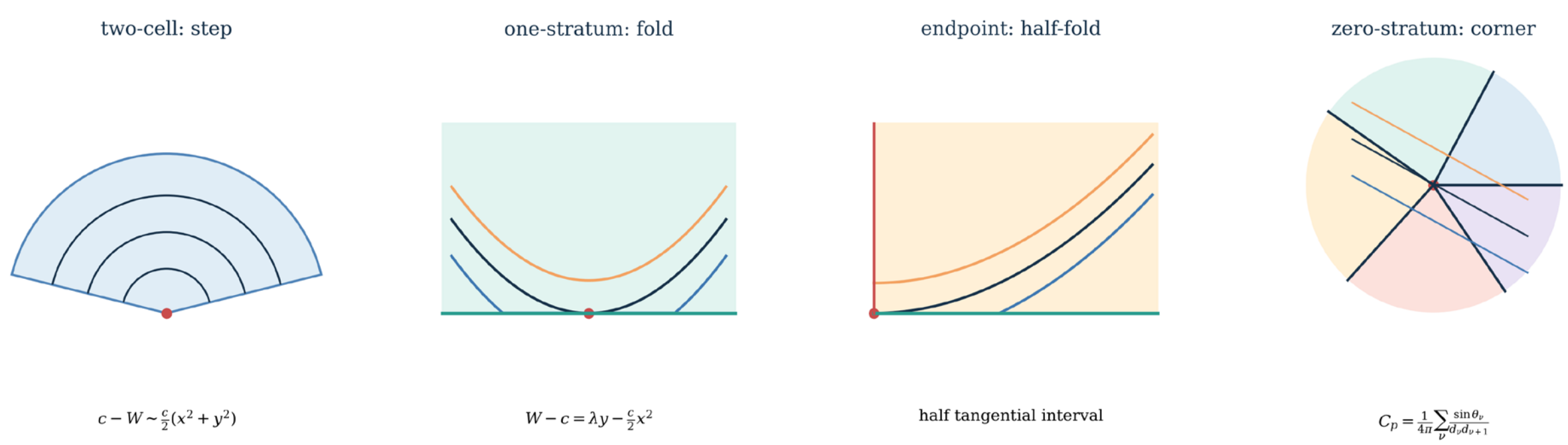

*Figure 3. The common local models for a two-cell step, an interior fold, an endpoint half-fold, and a zero-stratum corner. The exact coefficients are derived in Subsections 8.5–8.7.*

### 8.5. Two-cells: the step coefficient

Suppose $p = p_C \in \overline{C}$. Then $G_C = cp$, $c = |G_C|$, and (8.5) gives the exact cell identity

$$W\big(u(x)\big) = \frac{c}{\sqrt{1 + |x|^2}}. \qquad (8.8)$$

Let $\Theta_C(p)$ be the Euclidean angle of the tangent cone of $\overline{C}$ at $p$, with $\Theta_C(p) = 2\pi$ when $p \in C$. From (8.8),

$$c - W\big(u(x)\big) = \frac{c}{2}|x|^2 + O(|x|^4).$$

Consequently

$$\text{area}\{u \in C : c - \delta < W(u) \le c\} = \frac{\Theta_C(p)}{c}\,\delta + O(\delta^2).$$

After division by $4\pi$, the cell contribution to the density has the downward jump

$$J_C^{\text{step}} = \frac{\Theta_C(p)}{4\pi c}. \qquad (8.9)$$

The angle is positive for every incident two-cell. Hence step contributions at a common value have one sign and cannot cancel. The antipodal cell has the same coefficient.

### 8.6. One-strata: the fold coefficient

Let $p = p_E \in E$ be the stationary point of an open edge in $g_j^{\perp}$. Write

$$\beta = |g_j|, \qquad q = \frac{g_j}{\beta}, \qquad \alpha = H_E \cdot q, \qquad c = c_E.$$

Choose $x$ tangent to $E$ and signed $y$ in direction $q$. The two adjacent cells have gradients $H_E + g_j$ for $y > 0$ and $H_E - g_j$ for $y < 0$. Since $p$ is stationary on $E$, the tangent component of $H_E$ is $\alpha q$. With inward coordinates $y_+ = y \ge 0$ and $y_- = -y \ge 0$, (8.6) becomes

$$W - c = \lambda_\pm y_\pm - \frac{c}{2}x^2 + O\big(|x|^3 + |x|y_\pm + y_\pm^2\big), \qquad \lambda_\pm = \beta \pm \alpha. \qquad (8.10)$$

The two transverse alternatives in (8.10) govern the local fold.

Along the edge, arclength may be chosen so that $W(x,0) = c\cos x$. Assume first $\lambda \neq 0$. Solving the level equation $W - c = t$ for $y$ shows that the admissible $x$-interval changes at

$$x = \pm\sqrt{-\frac{2t}{c}} + O(t), \qquad t \uparrow 0.$$

Integration in $y$ contributes the reciprocal transverse derivative. For either sign of $\lambda$, the nonanalytic part of the unnormalised half-cell pushforward is

$$-\frac{2}{\lambda}\sqrt{\frac{2}{c}}\sqrt{(-t)_+} + O\big(|t|^{3/2}\big).$$

For $\lambda > 0$ this is the small interval removed when $t < 0$; for $\lambda < 0$ it is the small interval inserted, and the same signed formula results. Adding the two adjacent cells and dividing by $4\pi$ gives the local fold germ at $p$:

$$f_p(c+t) = A_p(t) - \frac{1}{2\pi}\sqrt{\frac{2}{c}}\left(\frac{1}{\beta+\alpha} + \frac{1}{\beta-\alpha}\right)\sqrt{(-t)_+} + O\big(|t|^{3/2}\big), \qquad (8.11)$$

where $A_p$ is real-analytic at zero. The absence of an order-$|t|$ one-sided term follows from the exact even edge restriction $c\cos x$: the two moving endpoints are opposite, and the odd first variation of the transverse Jacobian cancels between them. The coefficient reduces to

$$-\frac{1}{\pi}\sqrt{\frac{2}{c}}\frac{\beta}{\beta^2-\alpha^2}, \qquad (8.12)$$

and is nonzero whenever $\lambda_+\lambda_- \neq 0$. Its sign is

$$-\operatorname{sgn}(\beta^2-\alpha^2).$$

Thus $|\alpha| < \beta$ is the maximum-type fold familiar from orthotopes, whereas $|\alpha| > \beta$ is a stratified saddle and reverses the signed square-root germ. The term *fold* includes both regimes.

If $\lambda_+ = 0$ or $\lambda_- = 0$, the corresponding adjacent cell has $G_C = cp$; hence it contributes the step (8.9). The other half-cell still retains the fold appearing in (8.11). Thus a boundary-stationary two-cell produces an exact step–fold superposition, with the step dominant.

The antipodal point has the same germ. Distinct edge pairs attaining the same numerical value contribute the algebraic sum of their coefficients (8.12). Unlike the step coefficients, these coefficients need not have a common sign.

### *8.7. Zero-strata: the cyclic corner coefficient*

Let $p$ be a zero-stratum and suppose first that no incident open edge is stationary at $p$. List the unit tangent rays of the incident arrangement in positive cyclic order,

$$r_1, r_2, \dots, r_m, r_{m+1} = r_1,$$

and let $\theta_\nu \in (0,\pi)$ be the angle from $r_\nu$ to $r_{\nu+1}$. Let $C_\nu$ be the intervening two-cell and $q_\nu = G_{C_\nu} - cp$. Continuity of $W$ across each ray makes

$$d_\nu = q_{\nu-1} \cdot r_\nu = q_\nu \cdot r_\nu$$

well-defined. The nonstationarity assumption is exactly $d_\nu \neq 0$ for every $\nu$.

On the sector generated by $r_\nu, r_{\nu+1}$, write $x = sr_\nu + tr_{\nu+1}$, $s,t \geq 0$. Its Jacobian is $\sin\theta_\nu$, and the first jet of $W - c$ is

$$d_\nu s + d_{\nu+1} t.$$

The derivative jump in the pushforward density of a linear form $As + Bt$ on a quadrant is $1/(AB)$. This follows directly in all sign cases. If $A, B$ have the same sign, the level segment born at the origin has length parameter $|\rho|/|AB|$; if their signs are opposite, the level line crosses the quadrant and the endpoint transferred from one axis to the other has the same signed increment $1/(AB)$. Multiplication by the sector Jacobian gives

$$\frac{\sin\theta_\nu}{d_\nu d_{\nu+1}}.$$

To see that no hidden cutoff term enters the sector calculation leading to (8.13), choose a smooth localizer equal to one near the origin. The difference of any two such localizers is supported where the level map is a submersion and hence has an analytic pushforward. In the remaining neighbourhood, (8.6), the implicit function theorem, and the nonzero boundary derivatives $d_\nu$ move every sector endpoint by $O(\rho^2)$ beyond its linear displacement. They therefore change the two one-sided derivatives equally and do not alter their jump.

After summing sectors and normalizing spherical area, the exact local corner coefficient is

$$C_p^{\text{corner}} = \frac{1}{4\pi} \sum_{\nu=1}^{m} \frac{\sin\theta_\nu}{d_\nu d_{\nu+1}}. \qquad (8.13)$$

Thus $p$ contributes a corner if and only if the cyclic sum in (8.13) is nonzero.

If $d_\nu = 0$, the incident ray $r_\nu$ is stationary at its endpoint. Let $g_j^\perp$ be the great circle containing that ray, orient $q = g_j/|g_j|$ toward one adjacent sector, put $\beta = |g_j|$, and let $\alpha = H_E \cdot q$ for the sign pattern on the open edge issuing along $r_\nu$. The tangent coordinate now satisfies $x \geq 0$, so the interval calculation of Subsection 8.6 has half its interior-edge length. The two adjacent sectors contribute

$$-\frac{1}{4\pi}\sqrt{\frac{2}{c}}\left(\frac{1}{\beta+\alpha} + \frac{1}{\beta-\alpha}\right)\sqrt{(c-\rho)_+}. \qquad (8.14)$$

This coefficient is nonzero when $\beta \pm \alpha \neq 0$. If one denominator vanishes, the corresponding incident cell is stationary and Subsection 8.5 supplies a step; the other half-sector retains its fold. Multiple stationary incident rays contribute the sum of their one-sided coefficients. These are exactly the stronger mechanisms predicted by the common normal form. After those sectors are removed, the residual corner coefficient is the sum in (8.13) over the sectors for which both boundary derivatives remain nonzero.

At a simple vertex exactly two generators, say $g, h$, are active. The linear map $x \mapsto (z_1, z_2) = (g \cdot x, h \cdot x)$ is nonsingular. Write $\ell_p \cdot x = az_1 + bz_2$, and put $D = \left|\det_{T_pS^2}(g, h)\right|$. In the four quadrants the first variation is

$$|z_1| + |z_2| + az_1 + bz_2,$$

and $dx = D^{-1}dz_1dz_2$. Applying the quadrant calculation gives

$$\begin{aligned} 4\pi C_p^{\text{corner}} &= \frac{1}{D}\Big(\frac{1}{(1+a)(1+b)} + \frac{1}{(1-a)(1+b)} \\ &\quad + \frac{1}{(1-a)(1-b)} + \frac{1}{(1+a)(1-b)}\Big) \\ &= \frac{4}{D(1-a^2)(1-b^2)}. \end{aligned} \tag{8.15}$$

The excluded values $a = \pm 1$ or $b = \pm 1$ are exactly stationary incident edges. Thus every simple zero-stratum has a nonzero corner coefficient unless a stronger fold or step mechanism is already present. Exact corner cancellation is a phenomenon of non-simple vertices.

The local coefficients also recover the companion-paper box formulas exactly. For an orthogonal box, a one-stratum associated with generator $a_ie_i$ has $\beta = a_i$, $\alpha = 0$, and stationary value

$$d_i = \left(\sum_{k \neq i} a_k^2\right)^{1/2}.$$

Summing (8.12) over its four attained stationary points gives

$$-\frac{4}{\pi a_i}\sqrt{\frac{2}{d_i}},$$

the square-root coefficient obtained there. At a simple box vertex, $a = b = 0$ and $D = a_ia_j$ in (8.15), so the antipodal pair gives the derivative jump

$$\frac{2}{\pi a_i a_j}.$$

The global density, fold, and corner reductions are therefore consistent at all three singularity levels, without numerical verification.

If repeated side lengths make several box strata attain the same critical value, the observed coefficient is their local algebraic superposition as derived in Subsection 8.6 and aggregated according to Subsection 8.9; equal contributions therefore appear with their full multiplicity.

### *8.8. The universal corner assertion is false*

The cancellation condition in (8.13) is not vacuous. Consider the reduced full-dimensional family

$$g_1 = (1,0,0), \quad g_2 = (0,1,0), \quad g_3 = (1,1,0), \quad g_4 = (2\sqrt{3}, 0, 1). \tag{8.16}$$

At $p = (0,0,1)$, the active generators are $g_1, g_2, g_3$, and $c = W(p) = 1$. Writing tangent coordinates as $(x, y)$, the first variation (8.7) is

$$H(x,y) = 2\sqrt{3}\,x + |x| + |y| + |x + y|. \qquad (8.17)$$

Put $A = 2\sqrt{3}$. The six cyclic unit rays and their directional slopes are

$$\begin{array}{cc} r_\nu & H(r_\nu) \\ (1,0) & A+2 \\ (0,1) & 2 \\ (-1,1)/\sqrt{2} & (2-A)/\sqrt{2} \\ (-1,0) & 2-A \\ (0,-1) & 2 \\ (1,-1)/\sqrt{2} & (A+2)/\sqrt{2}. \end{array} \qquad (8.18)$$

Thus the local model (8.17) has precisely the six nonzero slopes recorded in (8.18).

The successive angle sines are $1, 2^{-1/2}, 2^{-1/2}, 1, 2^{-1/2}, 2^{-1/2}$.

Figure 4. Exact corner cancellation at a non-simple zero-stratum

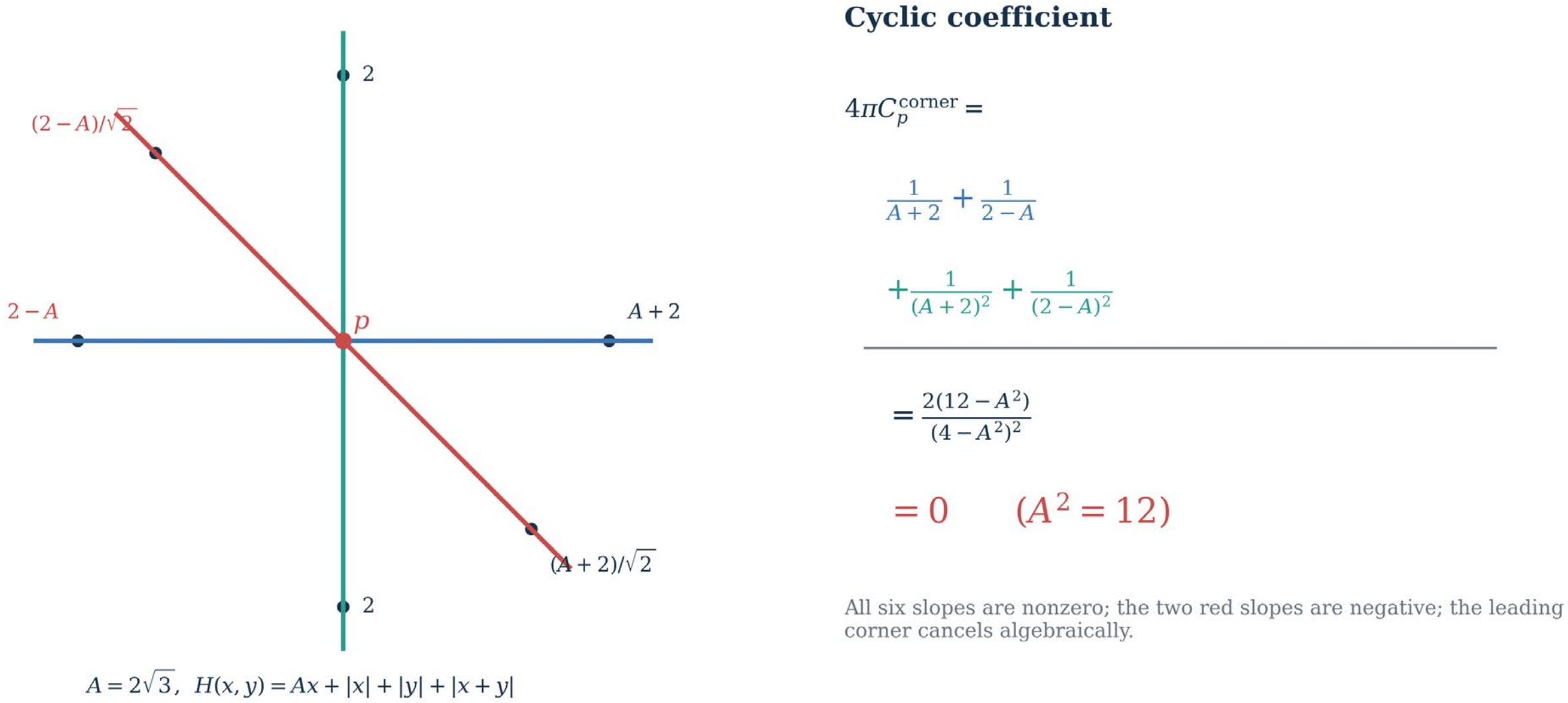

*Figure 4. The tangent arrangement, six nonzero directional slopes, and the cyclic coefficient for the family in (8.16). The following hand calculation, not the diagram, proves the cancellation.*

Substitution into the cyclic sum gives, by direct collection,

$$\begin{aligned} \sum_{\nu=1}^{6} \frac{\sin\theta_\nu}{d_\nu d_{\nu+1}} &= \frac{1}{A+2} + \frac{1}{2-A} + \frac{1}{(A+2)^2} + \frac{1}{(2-A)^2} \\ &= \frac{2(12 - A^2)}{(4 - A^2)^2} = 0. \end{aligned} \qquad (8.19)$$

No denominator vanishes because $A = 2\sqrt{3} \neq 2$. Hence no incident edge is stationary and no incident cell is stationary, but the zero-stratum corner coefficient is exactly zero. The family (8.16) is pairwise nonparallel and spans $\mathbb{R}^3$, so neither reduction nor degeneracy removes the example.

The cancelled value is isolated from every other mechanism. A direct application of (8.2)–(8.4), pairing antipodal strata, gives the exact zero-stratum values

$$\sqrt{\frac{2}{7}}, \qquad \frac{2}{\sqrt{13}}, \qquad 1, \qquad 2,$$

the attained open one-stratum values

$$\sqrt{\frac{56}{13}}, \qquad \sqrt{5}, \qquad \sqrt{9+4\sqrt{3}}, \qquad \sqrt{17+8\sqrt{3}},$$

and the admissible two-cell endpoint values

$$\sqrt{17+8\sqrt{3}}, \qquad \sqrt{21+8\sqrt{3}}.$$

The repeated value $\sqrt{17+8\sqrt{3}}$ is deliberate: it is simultaneously an attained one-stratum stationary value and an admissible two-cell endpoint, so its local density is governed by the step–fold superposition derived in Subsection 8.6 and aggregated according to Subsection 8.9.

Only the antipodal pair $\pm p$ attains the value 1. Hence no step, fold, or second corner contribution can mask the cancellation in (8.19).

This is a hand counterexample to the sentence "every zero-stratum contributes a corner." It also shows why a coefficient theorem is strictly stronger than an exponent table.

### *8.9. Superposition, dominance, and the generic theorem*

Choose a partition of unity subordinate to disjoint neighbourhoods of the finitely many stratified critical points and to their regular complement. The complement has an analytic pushforward by Subsection 8.3. At a fixed critical value $c$:

1. all step coefficients (8.9) are positive before the common downward sign, so a step can never cancel;
2. after the step germs are removed, all fold germs are given by (8.11) and their coefficients add algebraically;
3. after the fold germs are removed, all corner derivative jumps are given by (8.13) and add algebraically.

Therefore the dominance order

$$\text{step} \succ \text{square-root fold} \succ \text{corner} \qquad (8.20)$$

is universally correct as an order of local asymptotic scales. As a statement about the total density at a resonant value, the fold or corner level is visible exactly when its aggregate coefficient is nonzero. The qualification is necessary: step terms cannot cancel, but fold and corner terms can be signed.

For a generic reduced generator family, no three generators lie in a common plane through the origin, so every zero-stratum is simple. Formula (8.15) then makes its corner coefficient nonzero unless an incident edge is stationary. Formula (8.12) makes every isolated fold coefficient nonzero. The remaining exceptional conditions are equality in an incidence test or coincidence, beyond the forced antipodal duplication, of two values from the finite critical list. They are nontrivial real-analytic equations in the generators and are destroyed by an arbitrarily small independent perturbation. On each fixed combinatorial chart, their union is a finite union of proper real-analytic zero sets; it is closed relative to the chart, has empty interior, and has Lebesgue measure zero. Its complement is therefore open, dense, and of full measure. On that complement, every zero-stratum is a corner, every attained one-stratum stationary value is a square-root fold, and every admissible two-cell endpoint is a downward step.

### *8.10. Corrected reduced-zonotope singularity theorem*

**Theorem 8.1 (exact stratified width-density theorem).** For every full-dimensional reduced zonotope in $\mathbb{R}^3$, the density in (8.1) is real-analytic away from the finite set consisting of zero-stratum values, attained stationary values of open one-strata, including stratified saddles, and admissible stationary endpoints of two-cells. Its support is

$$\left[\min_{p \text{ zero-stratum}} W(p),\ \max_{C} |G_C|\right].$$

Incidence is decided exactly by (8.2), (8.3), and (8.4). The local step, interior-fold, endpoint-fold, and corner coefficients are respectively (8.9), (8.12), (8.14), and (8.13). Step contributions never cancel. Fold and corner contributions at a common value add algebraically, and the strongest nonzero aggregate term dominates in the order (8.20).

The unqualified form of Conjecture 7.13 of the companion paper [1] is false by (8.16)–(8.19). Its generic form is true, and its exact universal replacement is the incidence and coefficient theorem above. □

### *8.11. Geometric meaning and information boundary*

The theorem completes two lines of the paper. First, it extends the stratum principle of Section 6 beyond coordinate geometry. In both settings, lower-dimensional strata carry recoverable departures from a more universal aggregate, but the isolating operations differ: mixed-boundary Möbius inversion is specific to separable orthotopes, whereas the present section uses arrangement incidence and local co-area.

Second, it completes the passage initiated by Proposition 3.7 and Theorem 7.1. The Gram-matrix correlations that obstruct independent half-normal factorization do not obstruct the cellwise density calculus. Theorem 7.1 recovers the reduced generators from the ridges and gradient jumps of the direction-labelled function. The present theorem begins after that directional assignment has been forgotten and determines exactly how the scalar density is assembled from the same arrangement.

The resulting information boundary is precise:

$$\text{direction-labelled width function} \quad \Rightarrow \quad \text{reduced generators uniquely,}$$

whereas

$$\text{scalar width law} \quad \Rightarrow \quad \text{signed superposition of stratified contributions.}$$

Two-cell stationary values carry sign-rigid memory because their step coefficients cannot cancel. Isolated one-stratum stationary values have nonzero square-root coefficients, but resonant fold coefficients are signed. Simple zero-strata are visible unless a stronger incident mechanism occurs. At a non-simple zero-stratum, visibility is governed by the cyclic coefficient (8.13), which may vanish exactly.

Thus the section gives a complete forward calculus from a reduced generator arrangement to the singular structure of its width density. It does not prove that the singularity list reconstructs the entire arrangement, nor does the counterexample produce two distinct zonotopes with identical complete laws. Direction-labelled canonical uniqueness is proved in Section 7; complete scalar-law uniqueness for the full zonotopal or mixed linear–quadratic class is not asserted.

The correction of the former conjecture is substantive. The conjecture identified the generic exponent hierarchy but treated the presence of a zero-stratum as sufficient for a visible corner. The exact replacement proves universal sign rigidity for steps, universal local nonvanishing for isolated folds, universal nonvanishing at simple vertices, and the precise non-simple cancellation mechanism that remains possible. Density singularities are therefore geometric memory objects, but they are not automatic fingerprints: the observed law records their algebraic superposition.

### *8.12. Conjectural frontier and optimality problems*

The preceding theorems locate a remaining inverse problem that is both restricted enough to escape the counterexample of Subsection 7.1 and broad enough to contain all principal geometric classes treated here. Fix $d \geq 2$, a generator bound $N$, and an admissible quadratic rank

$$r \in \{0,2,3,\dots,d\}.$$

Consider canonical pairs $(G, A)$ with at most $N$ reduced generators and $\mathrm{rank} A = r$.

For the phrase *real-analytic exceptional set*, fix the number of generators, choose an ordering and signs locally, and work in a parameter chart where the generators are nonzero and pairwise nonparallel and the matrix has rank $r$. The exceptional sets in different charts are required to be invariant under the finite sign and permutation actions.

**Conjecture 8.2 (generic finite scalar tomography).** There exists a finite integer $M = M(d, N, r)$ such that every admissible fixed-generator-number chart contains a proper real-analytic exceptional set with the following property. If two canonical pairs, each lying outside the exceptional set of its respective chart, have at most $N$ reduced generators, quadratic rank $r$, and satisfy

$$\int_{S^{d-1}} \mathcal{W}_{G,A}(u)^k \, d\sigma(u) = \int_{S^{d-1}} \mathcal{W}_{\tilde{G},\tilde{A}}(u)^k \, d\sigma(u), \qquad 1 \leq k \leq M,$$

then they have the same number of generators and are orthogonally equivalent: there exist $Q \in O(d)$, a permutation $\pi$, and signs $\varepsilon_j \in \{-1,1\}$ such that

$$\tilde{A} = QAQ^T, \qquad \tilde{g}_{\pi(j)} = \varepsilon_j Q g_j.$$

The conjecture is deliberately generic. Unlabelled width laws are not injective on unrestricted centrally symmetric convex bodies by Subsection 7.1, while Subsection 8.8 shows that even within reduced zonotopes a local singularity coefficient may cancel on an exceptional family. On the positive side, Theorem 3.5 proves finite scalar-moment reconstruction for orthotopes, Propositions 7.4 and 7.5 prove law-level reconstruction for cylinders and ellipsoids, Theorem 7.3 proves finite reconstruction when directional ridge moments are retained, and Theorem 8.1 gives a generically noncancelling scalar signature for reduced zonotopes. These results motivate the conjecture but do not prove it.

Five quantitative problems accompany the conjecture.

1. **Optimal ridge-moment order.** Let $K(d, N_0)$ be the least integer for which $\mathsf{H}_0, \dots, \mathsf{H}_{K(d,N_0)}$ determine every canonical pair having at most $N_0$ reduced generators. Is the universal bound $2N_0 - 1$ in Theorem 7.3 sharp? What is the least generic order?
2. **Optimal scalar-moment order.** If Conjecture 8.2 holds, determine the least $M(d, N, r)$, and decide whether it can grow linearly with the dimension of the corresponding geometric parameter space.
3. **Optimal mixed-boundary data.** The transform isolating a fixed codimension-$k$ stratum uses its $2^k$ associated boundary patterns. What is the smallest global family of separable boundary patterns that reconstructs every codimension-$k$ coefficient?
4. **Optimal smoothing threshold.** Does the first Riesz mean admit a third deterministic coefficient at the transformed edge scale, with a remainder rigorously smaller than that scale? More generally, determine the least number of Riesz integrations for which such a coefficient exists.
5. **Cancellation geometry.** For a fixed nonsimple tangent arrangement, is the vanishing locus $C_p^{\text{corner}} = 0$ generically a codimension-one real-analytic hypersurface in generator space?

These questions preserve the information distinctions proved in the paper. They ask when less data suffice; none assumes that an ordinary unlabelled spectrum already contains direction-labelled width information.

# 9. Consequences and proved information boundaries

The preceding sections establish a common principle across stochastic width laws, spectral populations, and stratified spherical pushforwards. Leading aggregation can be universal even when the underlying geometry is fully recoverable. The information reappears either in a lower-order correction, in an isolated coordinate stratum, or in the incidence and coefficient data of a singular pushforward. Universality and geometric memory are therefore not opposing conclusions. They occur at different levels of the same decomposition.

### *9.1. Universality and nested boundary memory*

For orthotopes, Theorem 5.2 identifies a universal high-energy partition law: after normalization, the leading modal distribution converges to the Dirichlet law of squared spherical coordinates and loses the aspect ratio. Corollary 5.3 simultaneously shows that the associated maximizing directions recover the spherical width law of the companion article. This correspondence uses labelled modal directions; it does not identify that width law with an ordinary unlabelled spectrum.

Theorem 6.2 locates the first departure from universality. The singular part of its face-scale correction is supported on the simplex faces. Its coordinate moments determine the diagonal entries of $\mathsf{F}(K_a)$ up to the explicit factor displayed in Subsection 1.2, while the separate sine–cosine argument in Subsection 6.2 supplies the vanishing off-diagonal entries. The resulting projector-gradient tensor is basis-independent.

The underlying normal tensor $\mathsf{N}(K_a) = \int_{\partial K_a} \nu \otimes \nu \, dS$ is a classical interfacial Minkowski tensor in the cited normalization. The new role here is spectral: the directional face coefficient yields $\mathsf{F}(K_a)$, and that tensor explicitly recovers the orthotope.

At the next raw scale, Theorem 6.4 proves an obstruction rather than a missing estimate. Edge-sized jumps prevent an uncancelled deterministic third coefficient for the total sharp cutoff. The geometry has not disappeared: Theorems 6.5 and 6.6 use exact mixed-boundary Möbius inversion to isolate each coordinate stratum. In intrinsic dimension at least two they retain the next coefficient with a remainder strictly below its scale; a one-dimensional isolated stratum carries only the leading term and the stated remainder. The resulting hierarchy is recursive:

universal bulk → face memory → raw edge obstruction → exact stratum isolation
→ nested boundary memory.

Thus raw aggregation may obstruct a coefficient without destroying the geometric information carried by the corresponding stratum. Once a stratum of intrinsic dimension at least two is isolated, the same two-term bulk–boundary mechanism restarts in its lower dimension.

### *9.2. Recovery depends on the retained data*

The paper proves several reconstruction statements, but their inputs are distinct observables. Their information content becomes equivalent only after the corresponding reconstruction maps have been applied. For orthotopes, the distinction is exact:

1. the first $d$ scalar width moments or the $d$ nonconstant algebraic heat invariants of the complete Dirichlet trace determine the unordered side multiset;
2. the face-memory tensor or the isolated codimension-two spectral coefficients determine the coordinate-labelled side vector.

Labelled endpoints and indexed auxiliary laws give further, more directly resolved routes. These implications do not mean that an individual eigenvalue contains the interior of an auxiliary width law. The spectral endpoint map retains only the upper endpoint of each indexed law, and the ordinary spectrum additionally forgets the mode labels.

Beyond orthotopes, the data boundary becomes decisive. The complete direction-labelled function

$$\mathcal{W}_{G,A}(u) = \sum_j |g_j \cdot u| + 2\sqrt{u^T A u}$$

determines every canonical linear–quadratic pair by Theorem 7.1. Under a fixed generator bound, Theorem 7.3 reduces this continuum to finitely many direction-sensitive ridge moments and one quadratic residual tensor moment. These are distributional tensor data supported on the ridge arrangement; they are neither ordinary scalar moments of $\mathcal{W}_{G,A}(U)$ nor ordinary spectral aggregates.

For distinguished finite-dimensional classes, less information suffices. The unlabelled width law reconstructs a right circular cylinder, and finitely many even scalar moments reconstruct an ellipsoid. No corresponding universal claim is available for the full mixed class. Subsection 7.1 proves that unlabelled width laws are not injective even among smooth strictly convex centrally symmetric bodies. Accordingly,

$$\text{direction-labelled width data} \Rightarrow \text{canonical mixed geometry},$$

whereas

$$\text{unlabelled scalar width law on an unrestricted class} \not\Rightarrow \text{the body}.$$

This distinction is the exact boundary of the identifiability theorems proved here.

### 9.3. Stratified singularities as geometric memory

Section 8 shows how much structure survives after direction labels have been forgotten for a reduced zonotope. The global density is assembled from a finite great-circle arrangement. Its two-cells, one-strata, and zero-strata produce explicit step, fold, endpoint-fold, and corner coefficients. Steps are sign-rigid. Isolated folds and simple corners are locally nonzero. At resonant values, however, signed fold and corner contributions superpose, and the non-simple counterexample proves that a corner coefficient can cancel exactly.

This is the incidence counterpart of the Boolean stratum calculus in Section 6. For separable orthotopes, boundary-condition differences isolate coordinate strata before their coefficients are read. For reduced zonotopes, the great-circle incidence relations identify the local contributors, and the observed scalar density records their algebraic sum. In both cases, lower-dimensional strata carry geometric memory. What differs is whether an external transform isolates them or the scalar observable aggregates them.

Theorem 8.1 therefore supplies a complete forward calculus from a reduced generator arrangement to the singular structure of its scalar width density. It does not turn that density into a universal fingerprint. In particular, the cancellation example corrects the former universal corner assertion without producing two noncongruent zonotopes with identical complete laws. The exact scalar inverse problem remains separate from the proved direction-labelled reconstruction theorem.

### 9.4. The remaining frontier

Conjecture 8.2 asks whether bounded canonical linear–quadratic classes are generically determined, up to their natural orthogonal symmetries, by finitely many ordinary scalar moments. It seeks recovery from weaker, unlabelled scalar data than the direction-sensitive theorem of Section 7, but the two statements have different scopes and neither implies the other. The accompanying optimality questions ask how much ridge data, boundary-pattern data, smoothing, and noncancellation are actually necessary; they are not consequences of the present results.

The proved framework can therefore be summarized without conjectural identification:

| leading aggregation governs universality;<br>retained strata and directional incidence govern recoverable geometry. |
|---|

For orthotopes this principle yields explicit spectral reconstruction. For canonical mixed bodies it yields exact direction-sensitive identifiability. For reduced zonotopes it yields the full stratified singularity calculus and the precise cancellation boundary. Together these results separate what is universal, what is recoverable, and what remains genuinely open.

---

**Acknowledgements**. The material in this paper has been in development for many years and appears in print only now. I thank my family for their patience over that time. I am grateful to the early readers who criticised the early drafts; I owe a lasting debt to my late mathematics teacher, Amos Matalon, who showed me and seeded the love for advanced mathematics at an early age, insisted that it be taken seriously, and set me on the path to research. Any remaining errors are my own.

**Use of AI tools**. All derivations in this paper were carried out by hand. A large language model was used as an auxiliary tool in two respects: to run independent numerical checks confirming the closed forms stated here, and to assist with prior-art and bibliographic searching. Figures 1–4 are deterministic plots, produced with Matplotlib from the formulas proved in this paper; the plotting script was drafted with the assistance of a large language model and checked by the author. No figure in this paper is AI-generated imagery. The model produced no text and no mathematical argument in the paper, and the author takes full responsibility for its content.

**Declarations**. Funding: This research received no external funding. Competing interests: The author declares none. Data availability: The script generating Figures 1–4 is archived at 10.5281/zenodo.21827694. It is used for plotting and numerical confirmation only; no formula in this paper is derived from it.

---

## References

1. O. Abas, "Width distributions for rectangular boxes," arXiv:2608.03305 [math.MG] (2026).
2. S. Akiyama and T. Kamae, "Width deviation of convex polygons," *Discrete & Computational Geometry* **71** (2024), 1403–1428.
3. P. R. S. Antunes and P. Freitas, "Optimal spectral rectangles and lattice ellipses," *Proceedings of the Royal Society A* **469** (2013), article 20120492.
4. J. Borcea and P. Brändén, "Pólya–Schur master theorems for circular domains and their boundaries," *Annals of Mathematics* **170** (2009), 465–492.
5. P. Brändén and D. G. Wagner, "A converse to the Grace–Walsh–Szegő theorem," *Mathematical Proceedings of the Cambridge Philosophical Society* **147** (2009), 447–453.
6. H. Federer, *Geometric Measure Theory*, Grundlehren der mathematischen Wissenschaften 153, Springer, 1969.
7. S. R. Finch, "Width distributions for convex regular polyhedra," arXiv:1110.0671v2 (2016).
8. R. L. Frank and S. Larson, "Riesz means asymptotics for Dirichlet and Neumann Laplacians on Lipschitz domains," *Inventiones Mathematicae* **241** (2025), 999–1079.
9. R. J. Gardner, *Geometric Tomography*, 2nd ed., Encyclopedia of Mathematics and its Applications 58, Cambridge University Press, 2006.
10. K. Gittins and S. Larson, "Asymptotic behaviour of cuboids optimising Laplacian eigenvalues," *Integral Equations and Operator Theory* **89** (2017), 607–629.
11. M. Goresky and R. MacPherson, *Stratified Morse Theory*, Ergebnisse der Mathematik und ihrer Grenzgebiete, 3rd series, vol. 14, Springer, 1988.
12. J. H. Grace, "The zeros of a polynomial," *Proceedings of the Cambridge Philosophical Society* **11** (1902), 352–357.

13. E. Hlawka, “Über Integrale auf konvexen Körpern I,” *Monatshefte für Mathematik* **54** (1950), 1–36.
14. V. Ja. Ivrii, “Second term of the spectral asymptotic expansion of the Laplace–Beltrami operator on manifolds with boundary,” *Functional Analysis and Its Applications* **14** (1980), 98–106.
15. Z. Kabluchko, A. E. Litvak and D. Zaporozhets, “Mean width of regular polytopes and expected maxima of correlated Gaussian variables,” *Journal of Mathematical Sciences* **225** (2017), 770–787.
16. J. Kim, V. Yaskin and A. Zvavitch, “Distribution functions of sections and projections of convex bodies,” *Journal of the London Mathematical Society* (2) **95** (2017), no. 1, 52–72.
17. M. A. Klatt, G. Last, K. Mecke, C. Redenbach, F. M. Schaller and G. E. Schröder-Turk, “Cell shape analysis of random tessellations based on Minkowski tensors,” in *Tensor Valuations and Their Applications in Stochastic Geometry and Imaging*, Lecture Notes in Mathematics 2177, Springer, 2017, 385–421.
18. J. Lagacé, “Eigenvalue optimisation on flat tori and lattice points in anisotropically expanding domains,” *Canadian Journal of Mathematics* **72** (2020), no. 4, 967–987.
19. F. Lindemann, “Über die Zahl $\pi$,” *Mathematische Annalen* **20** (1882), 213–225.
20. P. Orlik and H. Terao, *Arrangements of Hyperplanes*, Grundlehren der mathematischen Wissenschaften 300, Springer, 1992.
21. R. Schneider, *Convex Bodies: The Brunn–Minkowski Theory*, 2nd expanded ed., Encyclopedia of Mathematics and its Applications 151, Cambridge University Press, 2014.
22. E. Sultanow and A. Hatziiliou, “Hearing the sides: Recovering a planar rectangle from eigenvalues,” arXiv:2511.23047v2 (2025).
23. E. Sultanow, A. Hatziiliou, C. May and N. Kratzke, “Hearing the edges: Recovering a 3D rectangular box from Dirichlet eigenvalues,” *Axioms* **15** (2026), article 284.
24. G. Szegő, “Bemerkungen zu einem Satz von J. H. Grace über die Wurzeln algebraischer Gleichungen,” *Mathematische Zeitschrift* **13** (1922), 28–55.
25. M. Tyaglov, “Generalized Hurwitz polynomials,” arXiv:1005.3032 (2010).
26. M. van den Berg, D. Bucur and K. Gittins, “Maximising Neumann eigenvalues on rectangles,” *Bulletin of the London Mathematical Society* **48** (2016), 877–894.
27. M. van den Berg and K. Gittins, “Minimising Dirichlet eigenvalues on cuboids of unit measure,” *Mathematika* **63** (2017), 469–482.
28. J. L. Walsh, “On the location of the roots of certain types of polynomials,” *Transactions of the American Mathematical Society* **24** (1922), 163–180.
29. A. G. Walters, “The distribution of projected areas of fragments,” *Proceedings of the Cambridge Philosophical Society* **43** (1947), no. 3, 342–347.
30. H. Weyl, “Ueber die asymptotische Verteilung der Eigenwerte,” *Nachrichten von der Gesellschaft der Wissenschaften zu Göttingen, Mathematisch-Physikalische Klasse* (1911), 110–117.